\documentclass[12pt,a4paper]{article}
\usepackage[a4paper, tmargin=3.4cm, bmargin=3.4cm, lmargin=2cm, rmargin=2cm, textheight=24cm, textwidth=16cm]{geometry}
\usepackage{setspace}
\usepackage{amsfonts}
\usepackage{epsfig}
\usepackage{latexsym}
\usepackage{amsthm}
\usepackage{mathtools,trimclip,lipsum}
\usepackage{hyperref}
\usepackage{xcolor}
\hypersetup{
colorlinks,
linkcolor={red!50!black},
citecolor={blue!50!black},
urlcolor={blue!80!black}
}
\usepackage{amsmath}
\usepackage{amssymb}

\usepackage[mathscr]{euscript}
 \let\mathscr\relax
\usepackage[scr]{rsfso}

\usepackage{array}
\usepackage{enumerate}
\usepackage{graphicx}
\usepackage{booktabs}
\usepackage{mathtools}
\usepackage{amsmath,blkarray,booktabs}

\newtheorem{theorem}{Theorem}[section]
\newtheorem{corollary}{Corollary}[section]
\newtheorem{proposition}{Propostion}[section]
\newtheorem{lemma}{Lemma}[section]

\newtheorem{re}{Remark}[section]
\newtheorem{nota}{Notation}[section]

\providecommand{\keywords}[1]
{
\small
\textbf{Keywords:} #1
}
\begin{document}

\title{ A method to optimize antipodal coloring span of graphs and its application}

\author{
Kush Kumar and Pratima Panigrahi
\\ \small Department of Mathematics, Indian Institute of Technology Kharagpur, India\\ \small e-mail: kushsingh029@gmail.com, pratima@maths.iitkgp.ac.in}


\maketitle
\begin{abstract}
In this article, we study radio \(k\)-colorings of simple connected graphs \(G\) with diameter \(d\), where a radio \(k\)-coloring \(g\) assigns non-negative integers to \(V(G)\) (vertices of \(G\)) such that \(|g(u) - g(v)| \geq 1 + k - d(u, v)\) for any two vertices \(u, v\) with \(1 \leq k \leq d\). The span of a radio \(k\)-coloring \(g\), expressed by \(rc_k(g)\), is the maximum integer assigned by \(g\), and the radio \(k\)-chromatic number \(rc_k(G)\) is the minimum span among all radio \(k\)-colorings of \(G\). A coloring \(g\) is minimal if \(rc_k(g) = rc_k(G)\). When \(k = d-1\), this coloring is known as the antipodal coloring, and \(rc_{d-1}(G)\) referred to as the antipodal number, is denoted by \(ac(G)\). We derive a sufficient condition for an antipodal coloring to be minimal and apply this criterion to determine the antipodal number of the generalized Petersen graph \(GP(n,1)\) for all \(n\) except when \(n \equiv 2 \pmod{8}\), and for toroidal grids \(T_{r,s} = C_r \square C_s\) when \(rs\) is even. Additionally, we establish a lower bound for \(ac(T_{r,s})\) when \(rs\) is odd.
\end{abstract}


\textbf{Mathematical Classification Code :}{ 05C12, 05C15, 05C76}

\section{Introduction}
Graph coloring is an important technique to address the frequency assignment problem, specifically in situations involving radio transmitters. The principal aim of this problem is assigning frequencies to transmitters such that interference reduces with minimized span. In this article, we consider a particular kind of radio $k$-coloring problem for a (simple connected) graph, which is antipodal coloring. For a simple connected graph $G$ with diameter $d$, a radio \(k\)-coloring \(g\) assigns non-negative integers to \(V(G)\) (the vertices of \(G\)) such that \(|g(u) - g(v)| \geq 1 + k - d(u, v)\) for any two vertices \(u, v\) with \(1 \leq k \leq d\). The span of $g$, expressed by $rc_k(g)$, is the largest integer assigned by $g$. The radio $k$-chromatic number, denoted as $rc_k(G)$, is the minimum span among all radio $k$-colorings of $G$. Notably, the coloring $g$ yielding the radio k-chromatic number must employ the color $0$. 

Furthermore, it is worth mentioning that when $k=1$, the radio $k$-coloring problem aligns with the traditional proper vertex coloring of graphs. Specific names are assigned to radio $k$-colorings for particular values of $k$: for instance, if $k=d$, then it is termed as radio coloring; if $k=d-1$, then it is called an antipodal coloring. Consequently, $rc_d(G)$ and $rc_{d-1}(G)$ respectively are known as the radio number (also denoted as $rn(G)$) and antipodal number (also denoted as $ac(G)$). A radio $k$-coloring $g$ of a graph $G$ is considered as a minimal radio $k$-coloring if $rc_k(g)=rc_k(G)$. Chartrand et al. \cite{chartrand2001radio} established the concept of radio $k$-coloring of graphs. Throughout the paper, we refer $|g(u)-g(v)|\geq d-d(u,v)$ as antipodal condition, and we consider only simple connected graphs.\\

\begin{nota}\label{Notation 1} For a $n$-vertex graph $G$, and a radio $k$-coloring $g$ of $G$, let $v_1, v_2,\ldots, v_n$ be an ordering of $V(G)$ such that $g(v_j)\leq g(v_{j+1})$ for all $j=1,2,\ldots, n-1$. Then $\epsilon_{j}$, $j=2,3,\ldots, n$, is defined as $\epsilon_{j}=g(v_j)- g(v_{j-1})-(1+k-d(v_j,v_{j-1}))$. We note that $\epsilon_{j}$ are non-negative integers.
\end{nota}

 The lemma presented below gives information about the radio $k$-chromatic number of an arbitrary graph. 
 \vspace{0.1cm}
 
 \begin{lemma}\label{Lemma 2}\cite{kola2010improved} Let $G$ be a graph of order $n$ and $g$ be a radio $k$-coloring of $G$. Then
\begin{equation}\label{equation 1}
        rc_{k}(g)=(n-1)(k+1)-\sum_{j=2}^{n}d(v_j,v_{j-1})+\sum_{j=2}^{n}\epsilon_j
    \end{equation}
where the vertices $v_j$ are arranged in the order given in Notation \ref{Notation 1}.
 \end{lemma}
 
 \begin{re}\label{Remark 2.1.}
{\rm In equation (\ref{equation 1}), for a given $k$ and a graph $G$, the term $(k+1)(n-1)$ is constant. So a radio $k$-coloring $g$ of $G$ is minimal if    $\sum_{j=2}^nd(v_{j-1},v_j)$  and $\sum_{j=2}^n{\epsilon}_j$ respectively attain the maximum and minimum value simultaneously among all possible radio $k$-colorings  of $G$}.
 \end{re}
 
The Cartesian product of two graphs $G$ and $H$, represented by $G\square H$, is a graph with vertex set $V(G)\times V(H)$, and two vertices $(u_1,v_1)$ and $(u_2,v_2)$ are adjacent in $G\square H$ if and only if $u_1=u_2$ and $v_1v_2\in E(H)$, or $u_1u_2\in E(G)$ and $v_1=v_2$. The graph $K_2\square C_n$ is the popular generalized Petersen graph $GP(n,1)$. The graph $C_r\square C_s$ is called the toroidal grid, and is denoted by $T_{r,s}$. Kchikech et al. \cite{kchikech2005radio} established some upper bounds for $rc_k(G\square H)$ as well as lower and upper bounds for $rc_k(P_n \square P_n)$ when $k$ is greater than or equal to $2n-3$. Chartrand et al. \cite{chartrand2002radio} derived the antipodal number for specific classes of paths and established some bounds for $ac(G)$ for any arbitrary graph $G$. The conjecture about the antipodal number of paths is given in \cite{chartrand2002radio}. Later, Khennoufa and Togni \cite{khennoufa2005note} have proved this conjecture. Also, Khennoufa and Togni \cite{khennoufa2011radio} determined $ac(Q_n)$ for $n$-dimensional hypercubes. In \cite{juan2012antipodal,chartrand2000radio}, the antipodal number of cycles is determined, while Saha and Panigrahi \cite{saha2012antipodal} calculated the antipodal number for powers of cycles. Niranjan et al. \cite{niranjan2020tree} provided an exact value for a specific class of trees $T$ and established bounds for $rc_k(T)$. In \cite{Basunia2021tree}, the authors determined the antipodal number of $m$-ary tree (for $m\geq3$ and height $h\geq3$), and constructed an optimal antipodal labeling. Gomathi et al. \cite{gomathi2022radio} obtained the bound for the antipodal number of honeycomb-derived networks. The radio $k$-chromatic number for corona of a graph with $K_1$ in \cite{niranjan2022corona}. Kim et al. \cite{kim2015radio} determined the radio number for the Cartesian product of a complete graph with a path, while Kola and Panigrahi \cite{kola2011radio} computed the radio number of $GP(n, 1)$. Kola and Panigrahi \cite{kola2015lower} provided a lower bound for $rc_k(G)$ for any arbitrary graph. They also verified that, for specific values of $r$, this lower bound coincides with $rn(C_r \square P_2)$. In \cite{Das2017}, Das et al. proposed a lower bound technique for radio $k$-coloring. The radio number of the toroidal grid $C_s \square C_s$ is determined by Morris-Rivera et al. \cite{morris2015radio}, and Saha and Panigrahi \cite{saha2013radio} determined $rn(C_s \square C_r)$ when $sr$ is even.
   
In the literature, we found no sufficient condition for the coloring of a graph to be minimal radio $k$-coloring. In Section 2, we give a sufficient condition for an antipodal coloring to be minimal. Applying this result, in Section 3, we determine the antipodal number of the generalized Petersen graph $GP(n,1)$ for all values of $n$ except the case when $n=4t+2$ with $t$ even. Moreover, for this case, we give an upper bound of $ac(GP(n,1))$. Finally, in Section 4, we obtain the antipodal number of toroidal grids $C_r\square C_s$ whenever $rs$ is even, and determine a lower bound for $ac(C_r\square C_s)$ whenever $rs$ is odd.
\section{\textbf{A Criterion for Minimal Antipodal Coloring}}
\label{criterion}
The fundamental result of this article is the subsequent theorem, which establishes a minimal antipodal coloring condition.

\begin{theorem}\label{Theorem 1.} Let $G$ be an $n$-vertex graph with diameter $d$. If an antipodal coloring $g$ of $G$ satisfies the condition $(a)$ or $(b)$ below according to $n$ is even or odd, then $g$ is a minimal antipodal coloring of $G$.
\begin{enumerate}
    \item [$(a)$] If $n$ is an even integer, then $g$ is such a mapping that  $d(v_j,v_{j+1})=d$, $d(v_{j+1},v_{j+2})=d(v_j,v_{j+2})+\epsilon_{j+1}+\epsilon_{j+2}$,  $d(v_{n-1},v_n)=d$ and $\epsilon_n=0$, for all odd $j$,  $1\leq j \leq n-3$.
    \item [$(b)$] If $n$ is an odd integer, then $g$ is such a mapping that  $d(v_j,v_{j+1})=d$, $d(v_{j+1},v_{j+2})=d(v_j,v_{j+2})+\epsilon_{j+1}+\epsilon_{j+2}$ \,for all odd $j$, $1\leq j \leq n-2$.
    
  Where $v_j$'s and $\epsilon_j$'s are as given in Notation \ref{Notation 1}.  
    
\end{enumerate}
\end{theorem}
\noindent
\textbf{\textit{Proof.}} Since $g$ is an antipodal coloring  of $G$, for three vertices $v_j,v_{j+1}$ and $v_{j+2}$, $1\leq j \leq n-3$, we have
\begin{equation}\label{equation 2}
    g(v_{j+1})-g(v_j)=d-d(v_j,v_{j+1})+\epsilon_{j+1},
    \end{equation}
    \begin{equation}\label{equation 3}
      g(v_{j+2})-g(v_{j+1})=d-d(v_{j+1},v_{j+2})+\epsilon_{j+2}.
\end{equation}
Now adding equations $(\ref{equation 2})$ and $(\ref{equation 3})$, we have
 \begin{equation}\label{equation 4}
     g(v_{j+2})-g(v_j)=2d-(d(v_j,v_{j+1})+d(v_{j+1},v_{j+2}))+(\epsilon_{j+1}+\epsilon_{j+2}).
 \end{equation}
 As $g$ is an antipodal coloring
 \begin{equation}\label{equation 5}
     g(v_{j+2})-g(v_j)\geq d-d(v_j,v_{j+2}).
 \end{equation}
From  $(\ref{equation 4})$ and $(\ref{equation 5})$, we get
\begin{equation}\label{equation 6}
    2d-(d(v_j,v_{j+1})+d(v_{j+1},v_{j+2}))+(\epsilon_{j+1}+\epsilon_{j+2})\geq d-d(v_j,v_{j+2})
    \end{equation}
    \begin{equation}\label{equation 7}
    \implies d(v_j,v_{j+1})+d(v_{j+1},v_{j+2})-d(v_j,v_{j+2})\leq d+(\epsilon_{j+1}+\epsilon_{j+2}).
\end{equation}
From equation $(\ref{equation 7})$, we get the following two inequalities:
\begin{equation}\label{equation 8}
    d(v_j,v_{j+1})+d(v_{j+1},v_{j+2})\leq d+d(v_j,v_{j+2})+\epsilon_{j+1}+\epsilon_{j+2}
\end{equation}
and
\begin{equation}\label{equation 9}
   \epsilon_{j+1}+\epsilon_{j+2} \geq d\left(v_{j}, v_{j+1}\right)+d\left(v_{j+1}, v_{j+2}\right)-\left(d\left(v_{j}, v_{j+2}\right)+d\right).
\end{equation}
\noindent
We note that  equality holds in  (\ref{equation 8}) and (\ref{equation 9})  for $d\left(v_{j}, v_{j+1}\right)=d$ and $d\left(v_{j+1}, v_{j+2}\right)=$ $d\left(v_{j}, v_{j+2}\right)+\epsilon_{j+1}+\epsilon_{j+2}$.\\

\noindent
\textbf{Case I:}  $n$ is an even integer.\\
In this case, the distance sum and the epsilon sum can be expressed as  $$\sum_{j=2}^{n} d\left(v_{j}, v_{j-1}\right)=\sum_{j=1, j\, \text { odd }}^{n-3}\left[d\left(v_{j}, v_{j+1}\right)+d\left(v_{j+1}, v_{j+2}\right)\right]+d\left(v_{n}, v_{n-1}\right)$$ and $$\sum_{j=2}^{n} \epsilon_{j}=\sum_{j=1, j \,\text { odd }}^{n-3}\left(\epsilon_{j+1}+\epsilon_{j+2}\right)+\epsilon_{n}.$$ 
Now, if $(a)$ holds true then  equality holds in both the equations (\ref{equation 8}) and (\ref{equation 9}), and so  the distance sum $\sum_{j=1, j\, \text { odd }}^{n-3}\left[d\left(v_{j}, v_{j+1}\right)+d\left(v_{j+1}, v_{j+2}\right)\right]+d\left(v_{n}, v_{n-1}\right)$ attains its maximum value and the epsilon sum $\sum_{j=1, j \,\text { odd }}^{n-3}\left(\epsilon_{j+1}+\epsilon_{j+2}\right)+\epsilon_{n}$ attains its minimum value simultaneously. Therefore, by Remark \ref{Remark 2.1.}, $g$ is a minimal antipodal coloring of $G$.\\

\noindent
\textbf{Case II:} $n$ is an odd integer.\\
In this case, the distance sum and the epsilon sum can be expressed as  $$\sum_{j=2}^{n} d\left(v_{j}, v_{j-1}\right)=\sum_{j=1, j\, \text { odd }}^{n-2}\left[d\left(v_{j}, v_{j+1}\right)+d\left(v_{j+1}, v_{j+2}\right)\right]$$ and $$\sum_{j=2}^{n} \epsilon_{j}=\sum_{j=1, j \,\text { odd }}^{n-2}\left(\epsilon_{j+1}+\epsilon_{j+2}\right).$$
Now, if $(b)$ holds true then equality holds in both (\ref{equation 8}) and (\ref{equation 9}), and the distance sum and epsilon sum respectively attain their maximum and minimum value simultaneously. Therefore, by Remark \ref{Remark 2.1.}, $g$ is a minimal antipodal coloring of $G$.\qed

\section{\textbf{Antipodal number of $GP(n,1)$}}
\label{ac(GP(n,1))}

In any graph $G$, two vertices $u$ and $v$ are said to be antipodal if the distance between them is equal to the diameter of $G$. We recall that $GP(n,1)$ consists of two cycles, $C_1:y_0y_1\ldots y_{n-1}y_0$ (say the inner cycle) and $C_2:x_0x_1\ldots x_{n-1}x_0$ (say the outer cycle) together with the edges $\{x_j,y_j\}$ for all $j=0,1,2,\ldots,n-1$. In $GP(n,1)$, each vertex has only one antipodal vertex if $n$ is even, and exactly two antipodal vertices if $n$ is odd.

The following lemma establish some fundamental properties of $GP(n,1)$.

\begin{lemma}\label{Lemma 2.2.}\cite{kola2011radio}  Let $v$ and $v^{\prime}$ be  antipodal vertices in $GP(n,1)$, where $n$ is even. Then for each vertex $w$ in $GP(n,1)$, $d(v,v^{\prime})=d(v,w)+d(w,v^{\prime})$ .
\end{lemma}

\begin{lemma}\label{Lemma 2.3.}\cite{kola2011radio} For an odd integer $n$, let $v^{\prime}$ and $v^{\prime\prime}$ be the two antipodal vertices of a vertex $u$ in $GP(n,1)$. Then for each vertex $w$ in $GP(n,1)$, we have either $d(u,v^{\prime})=d(u,w)+d(w,v^{\prime})$ or $d(u,v^{\prime\prime})=d(u,w)+d(w,v^{\prime\prime})$.
\end{lemma}

\begin{lemma}\label{Lemma 2.4.}\cite{kola2011radio}
For $GP(n,1)$, the diameter is given by
$$diam(\operatorname{GP}(n, 1))=\left\{\begin{array}{ll}\frac{n+2}{2} & \text { if } n \text { is even, } \\ \frac{n+1}{2} & \text { if } n \text { is odd. }\end{array}\right.$$
\end{lemma}
Since for any two vertices $u$ and $v$ in a graph $G$, $d(u,v)\geq 0$, we have the following remark:
\begin{re}\label{Remark Antipodal condition}
Let $G$ be a graph with diameter $d$ and $g$ be a vertex coloring of $G$. For $u,v\in V(G)$, if $g(u)-g(v)\geq d$, then $g$ satisfies the antipodal coloring condition for $u$ and $v$.
\end{re}
In the theorem below, we present the main result of this section.

\begin{theorem}\label{Theorem 2.} The antipodal number $ac(G P(n, 1))$ is 
$$
 \,\left\{\begin{array}{ll}
=\frac{n^{2}+3n-4}{4} & \text { if } n \equiv 0\pmod 4, \\
=\frac{n^{2}+2n-3}{4} & \text { if } n \equiv 1\pmod 4, \\
=\frac{n^{2}+5n-6}{4} & \text { if } n=4t+2, t \text { odd, } \\
\leq\frac{n^{2}+5n-6}{4} & \text { if } n=4t+2, t \text { even, }\\
=\frac{n^{2}-1}{4} & \text { if } n \equiv 3\pmod 4 .
\end{array}\right.
$$

\end{theorem}

\noindent
\textbf{\textit{Proof.}} Depending on the values of $n$, we consider four cases below, and where subscripts of $V\left(GP(n,1)\right)$ are taken modulo $n$.\\

\noindent
\textbf{Case I:} $n \equiv 0\pmod 4$.\\
First, we define vertex ordering of $GP(n,1)$ as 
$$
 v_{2j+1}=\left\{\begin{array}{ll}
{x}_{\frac{jn}{4}+1} & \text{ for } j=0,1,2,3,\\
    {y}_{\frac{n}{2}+\frac{(j-4)n}{4}+2} & \text{ for } j=4,5,6,7,\\
    {x}_{\frac{(j-8)n}{4}+3} & \text{ for } j=8,9,10,11,\\
    {y}_{\frac{n}{2}+\frac{(j-12)n}{4}+4} & \text{ for } j=12,13,14,15,\\
    \vdots  \\
    {y}_{\frac{n}{2}+\frac{(j-(n-4))n}{4}+\frac{n}{4}} & \text{ for } j=n-4,n-3,n-2,n-1,
\end{array}\right.
$$
and
$$
 v_{2j+2}=\left\{\begin{array}{ll}
{y}_{\frac{n}{2}+\frac{jn}{4}+1} & \text{ for } j=0,1,2,3,\\
    {x}_{\frac{(j-4)n}{4}+2} & \text{ for } j=4,5,6,7,\\
    {y}_{\frac{n}{2}+\frac{(j-8)n}{4}+3} & \text{ for } j=8,9,10,11,\\
    {x}_{\frac{(j-12)n}{4}+4} & \text{ for } j=12,13,14,15,\\
    \vdots  \\
    {x}_{\frac{(j-(n-4))n}{4}+\frac{n}{4}} & \text{ for } j=n-4,n-3,n-2,n-1.
\end{array}\right.
$$
From the above vertex ordering of $GP(n,1)$ and by Lemma \ref{Lemma 2.2.}, it can be  verified that $d(v_{2j+1},v_{2j+2})=\frac{n}{2}+1$ for $j\in \{0,1,2,\ldots,n-1\}$ and $d(v_{2j+2},v_{2j+3})=\frac{n}{4}+1$ for $j\in \{0,1,2,\ldots,n-2\}$. Then we get the vertex ordering of $GP(n,1)$ as $v_{1}, v_{2}, \ldots, v_{2 n}$ such that $\{d\left(v_{j}, v_{j+1}\right)\}_{j=1}^{2n-1}$ is an alternating sequence of $\frac{n}{2}+1$ and $\frac{n}{4}+1.$ Now, we give a coloring $g$ to $V\left(GP(n,1)\right)$ as: $g(v_1)=0$; $g(v_{j+1})=g(v_j)\,\, \text{ for } \,\,j=1,3,5,\ldots, 2n-1$; $g(v_{j+1})=g(v_j)+\frac{n}{4}+1 \,\, \text{ for }\,\, j=2,4,6,\ldots, 2n-2$.

We prove that $g$ is an antipodal coloring of $GP(n,1)$. Before this, we note that $v_j$ and $v_{j+1}$ are antipodal vertices for $j=1,3,\ldots,2n-3$, and by Lemma \ref{Lemma 2.2.}, $d(v_{j+1},v_{j})=d(v_{j+1},v_{j+2})+d(v_{j+2},v_{j})$. Then $d(v_{j+2},v_{j})=\left(\frac{n}{2}+1\right)-\left(\frac{n}{4}+1\right)=\frac{n}{4}$ for $j=1,3,\ldots,2n-3$. Also, by the ordering of $V\left(GP(n,1)\right)$, $d(v_{j+2},v_{j})\geq\frac{n}{4}$ for $j=2,4,\ldots,2n-2$.  For $j=1,2,\ldots,2n-2$, $g(v_{j+2})-g(v_j)=g(v_{j+2})-g(v_{j+1})+g(v_{j+1})-g(v_j)=\frac{n}{4}+1=1+\frac{n}{2}-\frac{n}{4}\geq1+(d-1)-d(v_{j+2},v_{j})$.

\noindent
For $j$ odd, since $g(v_{j+1})=g(v_j)$ and $g(v_{j+3})=g(v_{j+2})$, we get, 
    $g(v_{j+3})-g(v_j)=g(v_{j+3})-g(v_{j+2})+g(v_{j+2})-g(v_{j+1})+g(v_{j+1})-g(v_j)
    =g(v_{j+2})-g(v_{j+1})
    =\frac{n}{4}+1=1+\frac{n}{2}-\frac{n}{4}\geq 1+(d-1)-d(v_{j+3},v_{j})$
as $d(v_{j+3},v_{j})\geq d(v_{j+1},v_{j})-d(v_{j+3},v_{j+1})=\frac{n}{2}+1-\frac{n}{4}=\frac{n}{4}+1.$ For $j$ even,  $g(v_{j+2})=g(v_{j+1})$, and so we get,
 $g(v_{j+3})-g(v_j)=g(v_{j+3})-g(v_{j+2})+g(v_{j+2})-g(v_{j+1})+g(v_{j+1})-g(v_j)=g(v_{j+3})-g(v_{j+2})+g(v_{j+1})-g(v_j)=\frac{n}{4}+1+\frac{n}{4}+1 \geq 1+\frac{n}{2}-1\geq 1+(d-1)-d(v_{j+3},v_{j})$
as $d(v_{j+3},v_{j})\geq 1.$

\noindent
By the definition of $g$ and the ordering of $GP(n,1)$, the pair of vertices $v_j$ and $v_{j+1}$ safisfy the antipodal condition. Also, $g(v_{j+4})-g(v_j)=\sum_{l=1}^4[g(v_{j+l})-g(v_{j+l-1})]=2\left(\frac{n}{4}+1\right)=\frac{n}{2}+2\geq 1+(d-1)-d(v_{j+4},v_j)$. We also note that $g(v_{j+l})-g(v_{j})\geq d$ for $l\geq 5$, and so the antipodal condition holds true by Remark \ref{Remark Antipodal condition} for the pair $v_{j+l}$ and $v_{j}$ for $l\geq 5$.
Now the antipodal condition holds true between the pair of vertices $v_j$ and $v_{j+l}$ for all $l=1,2,\ldots,2n-j$, and hence $g$ is an antipodal coloring of $GP(n,1)$. Now using Lemma \ref{Lemma 2}, the span of $g$ is given as 
    $g(v_{2n})=(2n-1)\left(\frac{n}{2}+1\right)-\frac{2n}{2}\left(\frac{n}{2}+1\right)-\left(\frac{2n}{2}-1\right)\left(\frac{n}{4}+1\right)+n-1=\frac{n^2+3n-4}{4}.$
    
As by definition of $g$, $\epsilon_j=1$ for all $j=3,5,\ldots,2n-1$, and $\epsilon_j=0$ for all $j=2,4.\ldots,2n$. Now $d(v_j,v_{j+1})=\frac{n}{2}+1=diam(GP(n,1))$ and $d(v_{j+1},v_{j+2})=\frac{n}{4}+1=d(v_j,v_{j+2})+\epsilon_{j+1}+\epsilon_{j+2}$ for all $j=1,3,\ldots,2n-3$ and $d(v_{2n-1},v_{2n})=\frac{n}{2}+1$. Therefore, by Theorem \ref{Theorem 1.}, $g$ is a minimal antipodal coloring. Hence, the result follows.

\begin{figure}
    \centering
    \includegraphics[scale = 0.12]{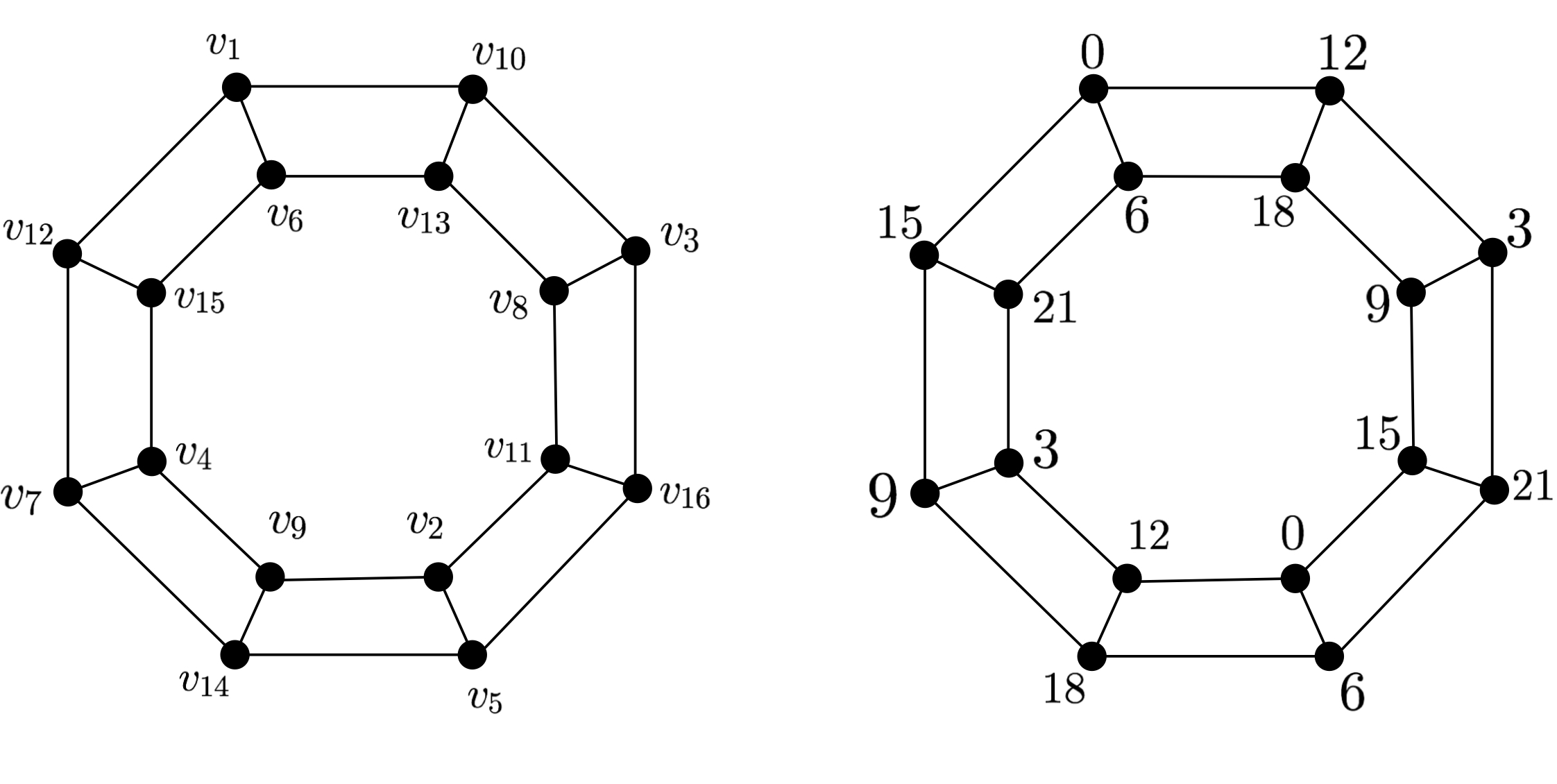}
    \caption{Ordering of vertices and mapping $g$ for $GP(8,1)$}
    \label{fig:1}
 \end{figure}

The ordering of vertices and mapping $g$ for $GP(8,1)$ are illustrated in Figure \ref{fig:1}.
\\
 
\noindent
\textbf{Case II:} $n \equiv 1\pmod 4$.\\
In this case, we define vertex ordering of $GP(n,1)$ as $v_{2j+1}={x}_{j\left(\frac{n-1}{4}\right)+1}$ and $v_{2j+2}={y}_{\left(\frac{n-1}{2}\right)+j\left(\frac{n-1}{4}\right)+1}$ for $j\in \{0,1,2,\ldots,n-1\}$. Then $d(v_j,v_{j-2})=\frac{n-1}{4}$ for all $j=3,4,\ldots,2n$. As $\frac{n-1}{4}$ and $n$ are co-prime, the sequence $v_1,v_3,\ldots,v_{2n-1}$ covers all the vertices of cycle $C_2$ and the sequence $v_2,v_4,\ldots,v_{2n}$ covers all the vertices of cycle $C_1$. Also, using the Lemma \ref{Lemma 2.3.}, we get $v_l$ and $v_{l-1}$ are antipodal vertices of each other for all $l\in \{2,4,\ldots, 2n\}$, and $d\left(v_{l-1}, v_{l+1}\right)+d\left(v_{l+1}, v_{l}\right)=d\left(v_{l-1}, v_{l}\right)$ for all $j\in \{2,4,\ldots,2n-2\}$. Since $diam\left(GP(n,1)\right)$ is $\frac{n+1}{2}$, for odd $n$, the distance sequence $d(v_1,v_2), d(v_2,v_3) \ldots, d(v_{2n-1},v_{2n})$ is an alternate sequence of $\frac{n+1}{2}$ and $\frac{n+3}{4}$. Now, we give a coloring $g$ to $V\left(GP(n,1)\right)$ as:
$g\left(v_{1}\right) =0$; 
$g\left(v_{j+1}\right) =g\left(v_{j}\right)\text{ for }j=1,3,5, \ldots, 2 n-1$; $g\left(v_{j+1}\right) =g\left(v_{j}\right)+\frac{n+3}{4} \text{ for } j=2,4,6, \ldots, 2 n-2$. Similarly as Case I,
we can verify that $g$ is an antipodal coloring. Now using Lemma \ref{Lemma 2}, the span of $g$ is given as
   $g(v_{2n})=(2n-1)\left(\frac{n+1}{2}\right)-\frac{2n}{2}\left(\frac{n+1}{2}\right)-\left(\frac{2n}{2}-1\right)\left(\frac{n+3}{4}\right)+n-1=\frac{n^2+2n-3}{4}$.
   
As by definition of $g$, $\epsilon_j=1$ for all $j=3,5,\ldots,2n-1$, and $\epsilon_j=0$ for all $j=2,4,\ldots,2n$. Now $d(v_j,v_{j+1})=\frac{n+1}{2}$, which is equal to the diameter of $GP(n,1)$, and $d(v_{j+1},v_{j+2})=\frac{n+3}{4}=\frac{n-1}{4}+1=d(v_j,v_{j+2})+\epsilon_{j+1}+\epsilon_{j+2}$ for all $j=1,3,\ldots,2n-3$, and $d(v_{2n-1},v_{2n})=\frac{n+1}{2}$. Therefore, by Theorem \ref{Theorem 1.}, $g$ is a minimal antipodal coloring, and hence the result follows.\\

\noindent
\textbf{Case III:} $n \equiv 2\pmod 4$.\\
In this case, we divide the proof into two subcases.

\noindent
\textbf{Subcase I:} Consider $n=4t+2$ and $t$ is odd.\\
First, we note that $t$ and $4t+2$ are co-prime. So  $\frac{n-2}{4}$ and $n$ are co-prime. 
In this case, we define an ordering of $V\left(GP(n,1)\right)$ as $v_{2j+1}={x}_{j\left(\frac{n-2}{4}\right)+1}$ and $v_{2j+2}={y}_{\left(\frac{n-1}{2}\right)+j\left(\frac{n-2}{4}\right)+1}$ for $j\in \{0,1,2,\ldots,n-1\}$. Then $d(v_j,v_{j-2})=\frac{n-2}{4}$ for all $j=3,4,\ldots,2n$. As $\frac{n-2}{4}$ and $n$ are co-prime, the sequence $v_1,v_3,\ldots,v_{2n-1}$ covers all the vertices of cycle $C_2$ and the sequence $v_2,v_4,\ldots,v_{2n}$ covers all the vertices of cycle $C_1$. Also, using the Lemma \ref{Lemma 2.2.}, we get $v_l$ and $v_{l-1}$ are antipodal vertices of each other for all $l\in \{2,4,\ldots, 2n\}$, and $d\left(v_{l-1}, v_{l+1}\right)+d\left(v_{l+1}, v_{l}\right)=d\left(v_{l-1}, v_{l}\right)$ for all $j\in \{2,4,\ldots,2n-2\}$. Since $diam\left(GP(n,1)\right)$ is $\frac{n}{2}+1$ for even $n$, the distance sequence $d(v_1,v_2), d(v_2,v_3) \ldots, d(v_{2n-1},v_{2n})$ is an alternate sequence of $\frac{n}{2}+1$ and $\frac{n+2}{4}+1$. Now, we give a coloring $g$ to $V\left(GP(n,1)\right)$ as: $g\left(v_{1}\right) =0$; 
$g\left(v_{j+1}\right) =g\left(v_{j}\right)  \text{ for } j=1,3,5, \ldots, 2 n-1$;
$g\left(v_{j+1}\right) =g\left(v_{j}\right)+\frac{n+2}{4}+1\text{ for }j=2,4,6, \ldots, 2 n-2.$ Similarly as Case I, we can verify that $g$ is an antipodal coloring. Now using Lemma \ref{Lemma 2}, the span of $g$ is given as
   $g(v_{2n})=(2n-1)\left(\frac{n}{2}+1\right)-\frac{2n}{2}\left(\frac{n}{2}+1\right)-\left(\frac{2n}{2}-1\right)\left(\frac{n+2}{4}+1\right)+2(n-1)=\frac{n^2+5n-6}{4}$.

As by definition of $g$, $\epsilon_j=2$ for all $j=3,5,\ldots,2n-1$, and $\epsilon_j=0$ for all $j=2,4,\ldots,2n$. Now $d(v_j,v_{j+1})=\frac{n}{2}+1=diam(GP(n,1))$ and $d(v_{j+1},v_{j+2})=\frac{n+2}{4}+1=\frac{n-2}{4}+2=d(v_j,v_{j+2})+\epsilon_{j+1}+\epsilon_{j+2}$ for all $j=1,3,\ldots,2n-3$, and $d(v_{2n-1},v_{2n})=\frac{n}{2}+1$. Therefore, by Theorem \ref{Theorem 1.}, $g$ is a minimal antipodal coloring, and hence the result follows.

\noindent
\textbf{Subcase II:} Consider $n=4t+2$ and $t$ is even.\\
\noindent
First, we note that $t+1$ and $4t+2$ are co-prime. So  $\frac{n+2}{4}$ and $n$ are co-prime. 
In this case, we define an ordering of $V\left(GP(n,1)\right)$ as $v_{2j+1}={x}_{j\left(\frac{n+2}{4}\right)+1}$ and $v_{2j+2}={y}_{\left(\frac{n-1}{2}\right)+j\left(\frac{n+2}{4}\right)+1}$ for $j\in \{0,1,2,\ldots,n-1\}$. Then $d(v_j,v_{j-2})=\frac{n+2}{4}$ for all $j=3,4,\ldots,2n$. As $\frac{n+2}{4}$ and $n$ are co-prime, the sequence $v_1,v_3,\ldots,v_{2n-1}$ covers all the vertices of cycle $C_2$ and the sequence $v_2,v_4,\ldots,v_{2n}$ covers all the vertices of cycle $C_1$. Also, using Lemma \ref{Lemma 2.2.}, we get $v_l$ and $v_{l-1}$ are antipodal vertices of each other for all $l\in \{2,4,\ldots, 2n\}$, and $d\left(v_{l-1}, v_{l+1}\right)+d\left(v_{l+1}, v_{l}\right)=d\left(v_{l-1}, v_{l}\right)$ for all $j\in \{2,4,\ldots,2n-2\}$. Since $diam\left(GP(n,1)\right)$ is $\frac{n}{2}+1$ for even $n$, the distance sequence $d(v_1,v_2), d(v_2,v_3) \ldots, d(v_{2n-1},v_{2n})$ is an alternate sequence of $\frac{n}{2}+1$ and $\frac{n+2}{4}$. Now, we give a coloring $g$ to $V\left(GP(n,1)\right)$ as:
$
g\left(v_{1}\right) =0$; $g\left(v_{j+1}\right) =g\left(v_{j}\right) \text{ for } j=1,3,5, \ldots, 2 n-1$; $g\left(v_{j+1}\right) =g\left(v_{j}\right)+\frac{n+2}{4}\text{ for } j=2,4,6, \ldots, 2 n-2.$
Similarly as Case I, we can verify that $g$ is an antipodal coloring. Now using Lemma \ref{Lemma 2}, the span of $g$ is given as
$g(v_{2n})=(2n-1)\left(\frac{n}{2}+1\right)-\frac{2n}{2}\left(\frac{n}{2}+1\right)-\left(\frac{2n}{2}-1\right)\left(\frac{n+2}{4}\right)+(n-1)=\frac{n^2+5n-6}{4}$.

As by definition of $g$, $\epsilon_j=1$ for all $j=3,5,\ldots,2n-1$, and $\epsilon_j=0$ for all $j=2,4,\ldots,2n$. Now $d(v_j,v_{j+1})=\frac{n}{2}+1=diam(GP(n,1))$ but $d(v_{j+1},v_{j+2})=\frac{n+2}{4}\neq d(v_j,v_{j+2})+\epsilon_{j+1}+\epsilon_{j+2}$ for all $j=1,3,\ldots,2n-3$. So we cannot apply Theorem \ref{Theorem 1.}, but the above antipodal coloring gives the upper bound for $ac(GP(n,1))$. Hence, the result follows.\\

\noindent
\textbf{Case IV:} $n \equiv 3\pmod 4$.\\
First, we note that $t+1$ and $4t+3$ are co-prime. So  $\frac{n+1}{4}$ and $n$ are co-prime. 
In this case, we define an ordering of $V\left(GP(n,1)\right)$ as $v_{2j+1}={x}_{j\left(\frac{n+1}{4}\right)+1}$ and $v_{2j+2}={y}_{\left(\frac{n-1}{2}\right)+j\left(\frac{n+1}{4}\right)+1}$ for $j\in \{0,1,2,\ldots,n-1\}$. Then $d(v_j,v_{j-2})=\frac{n+1}{4}$ for all $j=3,4,\ldots,2n$. As $\frac{n+1}{4}$ and $n$ are co-prime, the sequence $v_1,v_3,\ldots,v_{2n-1}$ covers all the vertices of cycle $C_2$ and the sequence $v_2,v_4,\ldots,v_{2n}$ covers all the vertices of cycle $C_1$. Also, using Lemma \ref{Lemma 2.3.}, we get $v_l$ and $v_{l-1}$ are antipodal vertices of each other for all $l\in \{2,4,\ldots, 2n\}$, and $d\left(v_{l-1}, v_{l+1}\right)+d\left(v_{l+1}, v_{l}\right)=d\left(v_{l-1}, v_{l}\right)$ for all $j\in \{2,4,\ldots,2n-2\}$. Since $diam\left(GP(n,1)\right)$ is $\frac{n+1}{2}$ for odd $n$, the distance sequence $d(v_1,v_2), d(v_2,v_3) \ldots, d(v_{2n-1},v_{2n})$ is an alternate sequence of $\frac{n+1}{2}$ and $\frac{n+1}{4}$. Now, we give a coloring $g$ to $V\left(GP(n,1)\right)$ as:
$g\left(v_{1}\right) =0$; $g\left(v_{j+1}\right) =g\left(v_{j}\right)\text{ for }j=1,3,5, \ldots, 2 n-1$; $g\left(v_{j+1}\right) =g\left(v_{j}\right)+\frac{n+1}{4}\text{ for } j=2,4,6, \ldots, 2 n-2.$
Similarly as Case I, we can verify that $g$ is an antipodal coloring, and using Lemma \ref{Lemma 2}, the span of $(g)$ is given as $g(v_{2n})=(2n-1)\left(\frac{n+1}{2}\right)-\frac{2n}{2}\left(\frac{n+1}{2}\right)-\left(\frac{2n}{2}-1\right)\left(\frac{n+1}{4}\right)=\frac{n^2-1}{4}.$

As by definition of $g$, $\epsilon_j=0$ for all $j=2,3,\ldots,2n$. Now $d(v_j,v_{j+1})=\frac{n+1}{2}=diam(GP(n,1))$ and $d(v_{j+1},v_{j+2})=\frac{n+1}{4}=d(v_j,v_{j+2})+\epsilon_{j+1}+\epsilon_{j+2}$ for all $j=1,3,\ldots,2n-3$ and $d(v_{2n-1},v_{2n})=\frac{n+1}{2}$. Therefore, by Theorem \ref{Theorem 1.}, $g$ is a minimal antipodal coloring, and hence the result follows.\qed

\section{\textbf{Antipodal number of $T_{r,s}$}}
 The vertex set of $T_{r,s}$ is defined as $V(T_{r,s})=\{(i,j):\, 0\leq i\leq r-1 \, \text{and} \, 0\leq j\leq s-1 \}$. In this context, all calculations involving the first and second components of $V\left(T_{r,s}\right)$ are performed using modulo $r$ and $s$, respectively. 

Below we state some properties of cycles and the Cartesian product of graphs, which we will use in the main results of this section.

\begin{re}\label{remark 2.1.1}\cite{saha2013radio} Let $i,j$ be the vertices of cycle $C_r$, then \\
$(i)$ $d(i,j)=\min\{|i-j|,r-|i-j|\}$.\\
$(ii)$ For any vertex $k$ of cycle $C_r$, $d(i,j)=d\left((i+k)\pmod{r},(j+k)\pmod{r})\right)$.   
\end{re}

\begin{proposition}\label{proposition 2.1.1.1}\cite{saha2013radio} Let $G$ and $H$ be simple connected graphs. Then\\
$(i)$ If $(g_1,h_1),(g_2,h_2)\in G\square H$, then
    $d_{G\square H}((g_1,h_1),(g_2,h_2))=d_{G}(g_1,g_2)+d_{H}(h_1,h_2).$\\
$(ii)$ $\text{diam}(G\square H)=\text{diam}(G)+\text{diam}(H)$.
\end{proposition}

\begin{corollary}\label{Corollary 2.1.1.}\cite{saha2013radio} For any two vertices $(i_1,j_1)$ and $(i_2,j_2)$ of $T_{r,s}$, then\\
$(i)$ $d((i_1,j_1),(i_2,j_2))=\min\{|i_1-i_2|,r-|i_1-i_2|\}+\min\{|j_1-j_2|,s-|j_1-j_2|\}.$\\
$(ii)$ For $(a,b)\in T_{r,s}$, $d((i_1,j_1),(i_2,j_2))=d((i_1+a,j_1+b),(i_2+a,j_2+b))$.
\end{corollary}

\begin{lemma}\label{Lemma 4.1.1} \cite{saha2013radio}  
Let \( n \) and \( r \) be positive integers with \( n > r \), and let \( \gcd(n, r) = p \). Then\\  
$(i)$ The elements of \( S = \{ir \pmod n : i = 0, 1, \ldots, \frac{n}{p} - 1\} \) are distinct. \\ 
$(ii)$ For any non-negative integer \( b \), the elements of \( S_b = \{(ir + b) \pmod n : i = 0, 1, \ldots, \frac{n}{p} - 1\} \) are also distinct.  
\end{lemma}
 
 Since the diameter of cycle $C_n$ is $\left\lfloor\frac{n}{2}\right\rfloor$, the following lemma is immediate from Proposition \ref{proposition 2.1.1.1}.
 
 \begin{lemma}\label{Lemma 2.1.01}
  The diameter of $T_{r,s}$ is $\left\lfloor\frac{r}{2}\right\rfloor+\left\lfloor\frac{s}{2}\right\rfloor$.
 \end{lemma}

The following lemma gives an ordering of $V\left(T_{r,s}\right)$ for $r\equiv 0\pmod 4$ and $s\equiv 0\pmod 4$, which will be useful to give a minimal antipodal coloring for $T_{r,s}$.

\begin{lemma}\label{Lemma for ac(T_{r,s}) when r,s=0 (mod 4)}
Let $r$ and $s$ be integers such that $r\equiv 0\pmod 4$ and $s\equiv 0\pmod 4$. Then there exists an ordering $v_1,v_2,\ldots,v_{rs}$ of $V\left(T_{r,s}\right)$ with the following properties:
\begin{enumerate}
    \item [$(a)$] $d(v_{2j},v_{2j-1})=\frac{r}{2}+\frac{s}{2}=diam(T_{r,s})$\quad for $j=1,2,\ldots,\frac{rs}{2}$.
    \item [$(b)$] $d(v_{2j+1},v_{2j})=\begin{cases}\frac{r}{4}+\frac{s}{4} & \text{ if } j=1,3,\ldots,\frac{rs}{2}-1, \\
    \frac{r}{4}+\frac{s}{4}+1 & \text{ if } j=2,4,\ldots,\frac{rs}{2}-2.
    \end{cases}$
    \item [$(c)$] $d(v_j,v_{j-2})=\begin{cases} \frac{r}{4}+\frac{s}{4} & \text{ if } j=3,4,7,8,\ldots,rs-1,rs, \\
    \frac{r}{4}+\frac{s}{4}-1 & \text{ if } j=5,6,9,10,\ldots,rs-3,rs-2.   
    \end{cases}$
    \item [$(d)$] $d(v_j,v_{j-3})=\begin{cases} \frac{r}{4}+\frac{s}{4} & \text{ if } j=4,8,\ldots,rs, \\
    \frac{r}{4}+\frac{s}{4}+1 & \text{ if } j=6,10,\ldots,rs-2, \\
    \frac{s}{2}+1 & \text{ if } j=5,7,9,\ldots,rs-1.
    \end{cases}$
\end{enumerate}
\end{lemma}
\noindent
\textbf{\textit{Proof.}} For each $j\in \{0,1,\ldots, s-1\}$, let $V_j=\{(i,j) : 0\leq i\leq r-1\}$. Then clearly $\{V_0, V_1,\ldots, V_{s-1}\}$ is a partition of $V(T_{r,s})$. Now we consider an ordering of $V\left(T_{r,s}\right)$ as:
    \(
        v_{4i+1}=\left(\frac{i(r+2)}{2},0\right),\,
        v_{4i+2}=\left(\frac{i(r+2)}{2}+\frac{r}{2},\frac{s}{2}\right),\,
        v_{4i+3}=\left(\frac{i(r+2)}{2}+\frac{3r}{4},\frac{3s}{4}\right),\,
        v_{4i+4}=\left(\frac{i(r+2)}{2}+\frac{r}{4},\frac{s}{4}\right) 
        \)
        for $0\leq i\leq r-1$, and 
       \( v_{4rj+l+1}=v_{l+1}+(j,j)\text{ for } 0\leq l\leq 4r-1 \text{ and }1\leq j\leq \frac{s}{4}-1.\)

    \noindent
    Now we show that the above ordering of $V\left(T_{r,s}\right)$ are distinct, and they covers all vertices of $T_{r,s}$. As $\frac{r}{2}$ is an even integer, so $\gcd\left(\frac{r}{2}+1,r\right)=1$. Therefore, by Lemma \ref{Lemma 4.1.1}, the elements of the set $\{v_{4i+1} : 0\leq i\leq r-1\}$ are distinct. Furthermore, applying Lemma \ref{Lemma 4.1.1}, it can be deduced that for every $j\in \{2,3,4\}$, all the elements in the set $\{v_{4i+j} : 0\leq i\leq r-1\}$ are distinct. Also, for every $j\in\{1,2,3,4\}$, $\{v_{4i+j} : 0\leq i\leq r-1\}=V_{\frac{ls}{4}}$ for some $l\in\{1,2,3,4\}$. So $\bigcup\limits_{j=1}^4\{v_{4i+j} : 0\leq i\leq r-1\}=\bigcup\limits_{j=1}^4V_{\frac{js}{4}}$. Hence, $v_{a}\neq v_{b}$ for any distinct integers $a,b\in\{1,2,\ldots,4r\}$. Further observe that $\{v_{4rj+l+1} : 0\leq l\leq 4r-1\}=\{v_{l+1}+(j,j) : 0\leq l\leq 4r-1\}=V_j\bigcup V_{\frac{s}{4}+j}\bigcup V_{\frac{s}{2}+j}\bigcup V_{\frac{3s}{4}+j}$ for each $j=0,1,\ldots \frac{s}{4}-1$. Therefore, the vertices in the ordering $v_1,v_2,\ldots,v_{rs}$ are distinct, and they cover all the vertices of $T_{r,s}$. Next, we have to show that the ordering of vertices defined above satisfies conditions $(a)-(d)$. We see that $v_{2j}=v_{2j-1}+\left(\frac{r}{2}+\frac{s}{2}\right)$ for all $j=1,2,\ldots,\frac{rs}{2}$. Then we get $d(v_{2j},v_{2j-1})=\frac{r}{2}+\frac{s}{2}$ for all $j=1,2,\ldots,\frac{rs}{2}$. Hence, condition $(a)$ is fulfilled. By applying Corollary \ref{Corollary 2.1.1.}, we compute the following distances:
    \begin{align*}
        d(v_{2j+1},v_{2j})=\begin{cases}\min \left\{\frac{3r}{4},r-\frac{3r}{4}\right\}+\min \left\{\frac{3s}{4},s-\frac{3s}{4}\right\}=\frac{r}{4}+\frac{s}{4} & \text{ if } j \text{ is odd},\\
            \min \left\{\frac{r}{4}+1,r-\frac{r}{4}-1\right\}+\min \left\{\frac{s}{4},s-\frac{s}{4}\right\}=\frac{r}{4}+\frac{s}{4}+1 & \text{ if } j \text{ is even}
        \end{cases}
    \end{align*} for $j \in \{1, 2, \ldots, 2r-1\}$. 
   Now, for $j\in\{3,4,\ldots,4r\}$,
    \begin{align*}
        d(v_{j},v_{j-2})=\begin{cases}\min \left\{\frac{3r}{4},r-\frac{3r}{4}\right\}+\min \left\{\frac{3s}{4},s-\frac{3s}{4}\right\}=\frac{r}{4}+\frac{s}{4} & \text{for } j=3,4,7,8,\ldots,4r-1,4r,\\
            \min \left\{\frac{r}{4}-1,r-\frac{r}{4}+1\right\}+\min \left\{\frac{3s}{4},s-\frac{3s}{4}\right\}=\frac{r}{4}+\frac{s}{4}-1 & \text{for } j=5,6,9,10,\ldots,4r-3,4r-2.
        \end{cases}
    \end{align*}
    For $j\in\{4,5,\ldots,4r\}$,
    \begin{align*}
        d(v_{j},v_{j-3})=\begin{cases}\min \left\{\frac{r}{4},r-\frac{r}{4}\right\}+\min \left\{\frac{s}{4},s-\frac{s}{4}\right\}=\frac{r}{4}+\frac{s}{4} & \text{ for } j=4,8,\ldots,4r,\\
        \min \left\{\frac{r}{4}+1,r-\frac{r}{4}-1\right\}+\min \left\{\frac{s}{4},s-\frac{s}{4}\right\}=\frac{r}{4}+\frac{s}{4}+1 & \text{ for } j=6,10,\ldots,4r-2,\\
            \min \left\{1,1-r\right\}+\min \left\{\frac{s}{2},s-\frac{s}{2}\right\}=\frac{s}{2}+1 & \text{ for } j=5,7,\ldots,4r-1.
        \end{cases}
    \end{align*}
    We have $v_{4r-2}=\left(r-1,\frac{s}{2}\right)$, $v_{4r-1}=\left(\frac{r}{4}-1,\frac{3s}{4}\right)$, $v_{4r}=\left(\frac{3r}{4}-1,\frac{s}{4}\right)$, $v_{4r+1}=\left(1,1\right)$, $v_{4r+2}=\left(\frac{r}{2}+1,\frac{s}{2}+1\right)$ and $v_{4r+3}=\left(\frac{3r}{4}+1,\frac{3s}{4}+1\right)$. Thus we get
    \begin{align*}
        &d(v_{4r-2},v_{4r+1})=\min \left\{r-2,r-r+2\right\}+\min \left\{\frac{s}{2}-1,s-\frac{s}{2}+1\right\}
        =\frac{s}{2}+1,\\
        &d(v_{4r-1},v_{4r+1})=\min \left\{\frac{r}{4}-2,r-\frac{r}{4}+2\right\}+\min \left\{\frac{3s}{4}-1,s-\frac{3s}{4}+1\right\}=\frac{r}{4}+\frac{s}{4}-1,\\
        &d(v_{4r},v_{4r+1})=\min \left\{\frac{3r}{4}-2,r-\frac{3r}{4}+2\right\}+\min \left\{\frac{s}{4}-1,s-\frac{s}{4}+1\right\}=\frac{r}{4}+\frac{s}{4}+1,\\
        &d(v_{4r-1},v_{4r+2})=\min \left\{\frac{r}{4}+2,r-\frac{r}{4}-2\right\}+\min \left\{\frac{s}{4}-1,s-\frac{s}{4}+1\right\}=\frac{r}{4}+\frac{s}{4}+1,\\
        &d(v_{4r},v_{4r+2})=\min \left\{\frac{r}{4}-2,r-\frac{r}{4}+2\right\}+\min \left\{\frac{s}{4}+1,s-\frac{s}{4}-1\right\}
        =\frac{r}{4}+\frac{s}{4}-1,\\
        &d(v_{4r},v_{4r+3})=\min \left\{2,r-2\right\}+\min \left\{\frac{s}{2}+1,s-\frac{s}{2}-1\right\}
        =\frac{s}{2}+1.
    \end{align*}
    Therefore condition $(b)$ holds for all $j\in\{1,2,\ldots,2r\}$, condition $(c)$ holds for $j\in\{3,4,\ldots,4r+1,4r+2\}$, and condition $(d)$ holds for $j\in\{4,5,\ldots,4r+1,4r+2,4r+3\}$. Since $v_{4rj+l+1}=v_{l+1}+(j,j)$ for $0\leq l\leq 4r-1$ and $1\leq j\leq \frac{s}{4}-1$, so using Corollary \ref{Corollary 2.1.1.}, conditions $(b)-(d)$ holds true for the rest of the values of $j$, and hence the result follows. \qed
\begin{theorem}\label{Theorem for ac(T_{r,s}) when r,s=0 (mod 4)} If $r\equiv 0\pmod 4$ and $s\equiv 0\pmod 4$, then
  $$ac(T_{r,s})=\frac{1}{8}\left(r^2s+rs^2+2rs-2r-2s-8\right).$$
\end{theorem}
\noindent
\textbf{\textit{Proof.}} Consider the ordering $v_1,v_2,\ldots,v_{rs}$ of $V\left(T_{r,s}\right)$ as given in the Lemma \ref{Lemma for ac(T_{r,s}) when r,s=0 (mod 4)}. We define $g : V(T_{r,s})\rightarrow \mathbb{N}\cup \{0\}$ such that 
\(
    \,g(v_1)=0;\,
    g(v_j)=g(v_{j-1}) \text{ for } j=2,4,6,\ldots, rs;\,
    g(v_j)=g(v_{j-1})+\left(\frac{r}{2}+\frac{s}{2}\right)-d(v_j,v_{j-1})\text{ for } j=3,7,\ldots, rs-1 \text{ and }\\
    g(v_j)=g(v_{j-1})+\left(\frac{r}{2}+\frac{s}{2}\right)-d(v_j,v_{j-1})+2 \text{ for } j=5,9,\ldots, rs-3.
\)
We show that $g$ is an antipodal coloring of $T_{r,s}$. From the definition of $g$, the antipodal condition holds true between the pair of vertices $v_{j}$ and $v_{j-1}$, for all $j=2,\ldots,rs$. We also note that $g(v_{j+l})-g(v_{j})\geq d$ for $l\geq 4$, and so the antipodal condition holds true by Remark \ref{Remark Antipodal condition} for the pair $v_{j+l}$ and $v_{j}$ for $l=4,5,\ldots,rs-j$ and $j=1,2,\ldots,rs-4$. Now, it remains to check that the antipodal condition satisfies between the pair of vertices $v_{j}$ and $v_{j-2}$ for $j=3,4,\ldots,rs$, and $v_{j}$ and $v_{j-3}$ for $j=4,5,\ldots,rs$. Depending on the value of $j$, we classify the analysis into the following four distinct cases:\\

\noindent
\textbf{Case I:} $j\equiv 0\pmod 4$.\\
From Lemma \ref{Lemma for ac(T_{r,s}) when r,s=0 (mod 4)}, $d(v_{j},v_{j-1})=\frac{r}{2}+\frac{s}{2}$, $d(v_{j-1},v_{j-2})=\frac{r}{4}+\frac{s}{4}$, $d(v_{j-2},v_{j-3})=\frac{r}{2}+\frac{s}{2}$, $d(v_{j},v_{j-2})=\frac{r}{4}+\frac{s}{4}$ and $d(v_{j},v_{j-3})=\frac{r}{4}+\frac{s}{4}$.\\ Therefore,
\(
    \,g(v_{j})-g(v_{j-2})=g(v_{j})-g(v_{j-1})+g(v_{j-1})-g(v_{j-2})
    =0+\left(\frac{r}{2}+\frac{s}{2}\right)-d(v_{j-1},v_{j-2})
    =\left(\frac{r}{2}+\frac{s}{2}\right)-\left(\frac{r}{4}+\frac{s}{4}\right)
    =diam(T_{r,s})-d(v_{j},v_{j-2}) \text{ and }
    g(v_{j})-g(v_{j-3})=g(v_{j})-g(v_{j-2})+g(v_{j-2})-g(v_{j-3})
    =\left(\frac{r}{4}+\frac{s}{4}\right)+0
    =\left(\frac{r}{2}+\frac{s}{2}\right)-\left(\frac{r}{4}+\frac{s}{4}\right)
    = diam(T_{r,s})-d(v_{j},v_{j-3}).
\)\\

\noindent
\textbf{Case II:} $j\equiv 1\pmod 4$.\\
From Lemma \ref{Lemma for ac(T_{r,s}) when r,s=0 (mod 4)}, $d(v_{j},v_{j-1})=\frac{r}{4}+\frac{s}{4}+1$, $d(v_{j-1},v_{j-2})=\frac{r}{2}+\frac{s}{2}$, $d(v_{j-2},v_{j-3})=\frac{r}{4}+\frac{s}{4}$, $d(v_{j},v_{j-2})=\frac{r}{4}+\frac{s}{4}-1$ and $d(v_{j},v_{j-3})=\frac{s}{2}+1$.\\ Therefore,
\(
    \, g(v_{j})-g(v_{j-2})=g(v_{j})-g(v_{j-1})+g(v_{j-1})-g(v_{j-2})
    =\left(\frac{r}{2}+\frac{s}{2}\right)-d(v_{j},v_{j-1})+2+0
    =\left(\frac{r}{2}+\frac{s}{2}\right)-\left(\frac{r}{4}+\frac{s}{4}+1\right)+2
    =\left(\frac{r}{2}+\frac{s}{2}\right)-\left(\frac{r}{4}+\frac{s}{4}-1\right)
    =diam(T_{r,s})-d(v_{j},v_{j-2}) \text{ and }
    g(v_{j})-g(v_{j-3})=g(v_{j})-g(v_{j-2})+g(v_{j-2})-g(v_{j-3})
    =\left(\frac{r}{4}+\frac{s}{4}+1\right)+\left(\frac{r}{2}+\frac{s}{2}\right)-d(v_{j-2},v_{j-3})
    =\left(\frac{r}{4}+\frac{s}{4}+1\right)+\left(\frac{r}{2}+\frac{s}{2}\right)-\left(\frac{r}{4}+\frac{s}{4}\right)
    =\frac{r}{2}+\frac{s}{2}+1
    \geq diam(T_{r,s})-d(v_{j},v_{j-3}).
\)\\

\noindent
\textbf{Case III:} $j\equiv 2\pmod 4$.\\
From Lemma \ref{Lemma for ac(T_{r,s}) when r,s=0 (mod 4)}, $d(v_{j},v_{j-1})=\frac{r}{2}+\frac{s}{2}$, $d(v_{j-1},v_{j-2})=\frac{r}{4}+\frac{s}{4}+1$, $d(v_{j-2},v_{j-3})=\frac{r}{2}+\frac{s}{2}$, $d(v_{j},v_{j-2})=\frac{r}{4}+\frac{s}{4}-1$ and $d(v_{j},v_{j-3})=\frac{r}{4}+\frac{s}{4}+1$.\\ Therefore,
\(
  \,  g(v_{j})-g(v_{j-2})=g(v_{j})-g(v_{j-1})+g(v_{j-1})-g(v_{j-2})=0+\left(\frac{r}{2}+\frac{s}{2}\right)-d(v_{j-1},v_{j-2})+2=\left(\frac{r}{2}+\frac{s}{2}\right)-\left(\frac{r}{4}+\frac{s}{4}+1\right)+2=\frac{r}{4}+\frac{s}{4}+1
    =diam(T_{r,s})-d(v_{j},v_{j-2}) \text{ and }
    g(v_{j})-g(v_{j-3})=g(v_{j})-g(v_{j-2})+g(v_{j-2})-g(v_{j-3})
    =\left(\frac{r}{4}+\frac{s}{4}+1\right)+0
    \geq \left(\frac{r}{2}+\frac{s}{2}\right)-\left(\frac{r}{4}+\frac{s}{4}+1\right)
    \geq diam(T_{r,s})-d(v_{j},v_{j-3}).
\)\\

\noindent
\textbf{Case IV:} $j\equiv 3\pmod 4$.\\
From Lemma \ref{Lemma for ac(T_{r,s}) when r,s=0 (mod 4)}, $d(v_{j},v_{j-1})=\frac{r}{4}+\frac{s}{4}$, $d(v_{j-1},v_{j-2})=\frac{r}{2}+\frac{s}{2}$, $d(v_{j-2},v_{j-3})=\frac{r}{4}+\frac{s}{4}+1$, $d(v_{j},v_{j-2})=\frac{r}{4}+\frac{s}{4}$ and $d(v_{j},v_{j-3})=\frac{s}{2}+1$.\\ Therefore,
\(
   \, g(v_{j})-g(v_{j-2})=g(v_{j})-g(v_{j-1})+g(v_{j-1})-g(v_{j-2})=\left(\frac{r}{2}+\frac{s}{2}\right)-d(v_{j},v_{j-1})+0=\left(\frac{r}{2}+\frac{s}{2}\right)-\left(\frac{r}{4}+\frac{s}{4}\right)
    =diam(T_{r,s})-d(v_{j},v_{j-2}) \text{ and }
    g(v_{j})-g(v_{j-3})=g(v_{j})-g(v_{j-2})+g(v_{j-2})-g(v_{j-3})
    =\left(\frac{r}{4}+\frac{s}{4}\right)+\left(\frac{r}{2}+\frac{s}{2}\right)-d(v_{j-2},v_{j-3})
    =\left(\frac{r}{4}+\frac{s}{4}\right)+\left(\frac{r}{2}+\frac{s}{2}\right)-\left(\frac{r}{4}+\frac{s}{4}+1\right)
    =\frac{r}{2}+\frac{s}{2}-1
    \geq diam(T_{r,s})-d(v_{j},v_{j-3}).
\)


\noindent
Hence, $g$ satisfies the antipodal condition between every pair of vertices $v_{j}$ and $v_{j-l}$ for $1\leq l\leq rs$ and $l+1\leq j\leq rs$. Therefore, $g$ is an antipodal coloring. As by definition of $g$, $\epsilon_{j}=0$ for all $j=2,3,4,6,7,8,10,\ldots,rs-4,rs-2,rs-1,rs$ and $\epsilon_{j}=2$ for all $j=5,9,\ldots,rs-3$. Now using Lemma \ref{Lemma 2}, the span of $g$ is given as $g(v_{rs})=(rs-1)\left(\frac{r}{2}+\frac{s}{2}\right)-\frac{rs}{2}\left(\frac{r}{2}+\frac{s}{2}\right)-\frac{rs}{4}\left(\frac{r}{4}+\frac{s}{4}\right)-\left(\frac{rs}{4}-1\right)\left(\frac{r}{4}+\frac{s}{4}+1\right)+2\left(\frac{rs}{4}-1\right) 
    =\frac{1}{8}\left(r^2s+rs^2+2rs-2r-2s-8\right).$

\noindent
Now $d(v_j,v_{j+1})=\frac{r}{2}+\frac{s}{2}=diam(T_{r,s})$ for all $j=1,3,5,\ldots,rs-3$. Also, $d(v_{j+1},v_{j+2})=\frac{r}{4}+\frac{s}{4}=d(v_j,v_{j+2})+\epsilon_{j+1}+\epsilon_{j+2}$ for all $j=1,5,9,\ldots,rs-3$, $d(v_{j+1},v_{j+2})=\frac{r}{4}+\frac{s}{4}+1=d(v_j,v_{j+2})+\epsilon_{j+1}+\epsilon_{j+2}$ for all $j=3,7,\ldots,rs-5$ and $d(v_{rs-1},v_{rs})=\frac{r}{2}+\frac{s}{2}=diam(T_{r,s})$. Therefore, by Theorem \ref{Theorem 1.}, $g$ is a minimal antipodal coloring, and hence the result follows.\qed
\begin{lemma} \label{Lemma for ac(T_{r,s}) when r=1,s=0 (mod 4) or r=3,s=2 (mod 4)}
Let $r$ and $s$ be integers such that $r\equiv 1\pmod 4$ and $s\equiv 0\pmod 4$, or $r\equiv 3\pmod 4$ and $s\equiv 2\pmod 4$. Then there exists an ordering $v_1,v_2,\ldots,v_{rs}$ of $V\left(T_{r,s}\right)$ with the following properties
\begin{enumerate}
    \item [$(a)$] $d(v_{2j},v_{2j-1})=\frac{r-1}{2}+\frac{s}{2}=diam(T_{r,s})$\quad for $j=1,2,\ldots,\frac{rs}{2}.$
    \item [$(b)$] $d(v_{2j+1},v_{2j})=\frac{r+s+3}{4}$\quad for  $j=1,2,\ldots,\frac{rs}{2}-1.$
    \item [$(c)$] $d(v_{j},v_{j-2})= \frac{r+s-1}{4}$\quad for $j=3,4,\ldots,rs.$ 
    \item [$(d)$] $d(v_j,v_{j-3})\begin{cases}= \frac{r+s-1}{4} & \text{ if } j=4,6,8,\ldots,rs, \\
    \geq1 & \text{ if } j=5,7,\ldots,rs-1.
    \end{cases}$
\end{enumerate}
\end{lemma}
\noindent
\textbf{\textit{Proof.}} \textbf{Case I:} $r\equiv 1\pmod 4$ and $s\equiv 0\pmod 4$.\\
Consider an ordering of $V\left(T_{r,s}\right)$ as:
    $
        v_{4i+1}=\left(\frac{i(r+1)}{2},0\right),\,
        v_{4i+2}=\left(\frac{i(r+1)}{2}+\frac{r-1}{2},\frac{s}{2}\right),\,\\
        v_{4i+3}=\left(\frac{i(r+1)}{2}+\frac{3r+1}{4},\frac{3s}{4}\right),\,
        v_{4i+4}=\left(\frac{i(r+1)}{2}+\frac{r-1}{4},\frac{s}{4}\right) \text{ for } 0\leq i\leq r-1$ and 
    $
        v_{4rj+l+1}=v_{l+1}+(j,j)\text{ for } 0\leq l\leq 4r-1 \text{ and }1\leq j\leq \frac{s}{4}-1.
    $
    Subsequently, we proceed similarly as in Lemma \ref{Lemma for ac(T_{r,s}) when r,s=0 (mod 4)}, and establish that the lemma is valid.\\

    \noindent
    \textbf{Case II:} $r\equiv 3\pmod 4$ and $s\equiv 2\pmod 4$.\\
    Consider an ordering of $V\left(T_{r,s}\right)$ as
    \begin{align*}
        v_{2rj+2i+1}&=\begin{cases}\left(\frac{i(r+1)}{4},\frac{j(s-2)}{4}\right) & \text{ if } i\equiv0\pmod2, \\
    \left(\frac{i(r+1)}{4},\frac{(j+1)(s-2)}{4}\right) & \text{ if } i\equiv1\pmod2
    \end{cases}\\
        v_{2rj+2i+2}&=\begin{cases}\left(\frac{i(r+1)}{4}+\frac{r+1}{2},\frac{j(s-2)}{4}+\frac{s}{2}\right) & \text{ if } i\equiv0\pmod2, \\
    \left(\frac{i(r+1)}{4}+\frac{r+1}{2},\frac{(j+1)(s-2)}{4}+\frac{s}{2}\right) & \text{ if } i\equiv1\pmod2
    \end{cases}
    \end{align*}
    for $i=0,1,\ldots,r-1$ and $j=0,1,\ldots,\frac{s}{2}-1$. Now we will show that the above ordering of $V\left(T_{r,s}\right)$ are distinct, and they cover all the vertices of $T_{r,s}$. If possible, assume $v_a=v_b$ for distinct integers $a$ and $b$ from the set $\{1,2,\ldots,rs\}$. If both $a$ and $b$ are even, then using division algorithm $a$ and $b$ can be written as $a=2rj_1+2i_1$ and $b=2rj_2+2i_2$ for integers $i_1, i_2, j_1$ and $j_2$ with $i_1, i_2 <r$. Now $v_a=v_b$ gives $\frac{(i_1-i_2)(r+1)}{4}\equiv0 \pmod{r}$ and $\frac{(j_1-j_2)(s-2)}{4}\equiv0 \pmod{s}$, whenever both $i_1$, $i_2$ are even integers. This implies that $i_1=i_2$ and $j_1=j_2$ as $\gcd\left(\frac{r+1}{4},r\right)=1$ and $\gcd\left(\frac{s-2}{4},s\right)=1 \text{ or }2$. Hence $v_a\neq v_b$ whenever $i_1$, $i_2$ are even integers. Similarly, $v_a\neq v_b$ whenever both $i_1$, $i_2$ are odd integers. Without losing generality, assume $i_1$ is an even and $i_2$ is an odd integer, so the first component of $v_a$ and $v_b$ are not same. Hence $v_a\neq v_b$ whenever both $a$ and $b$ are even integers. Also, from the ordering of vertices, $v_l=v_{l+1}-\left(\frac{r+1}{2},\frac{s}{2}\right)$ for every odd integer $l$. Hence $v_a\neq v_b$ whenever both $a$ and $b$ are odd integers. Then it remains to prove that $v_a$ and $v_b$ are distinct whenever $a$ is an even integer and $b$ is an odd integer. So assume $a=2rj_1+2i_1$ and $b=2rj_2+2i_2+1$, $i_1, i_2<r$. The second components of $v_a$ and $v_b$ give $\left(\frac{(j_2-j_1)(s-2)}{4}+\frac{s}{2}\right)\equiv0\pmod{s}$ whenever $i_1$ and $i_2$ are even integers, which is impossible as $\frac{(j_2-j_1)(s-2)}{4}+\frac{s}{2}$ is odd integer and $0\leq j_1, j_2\leq \frac{s}{2}-1$. Hence $v_a\neq v_b$ whenever $i_1$, $i_2$ are even integers. Similarly, $v_a\neq v_b$ whenever both $i_1$, $i_2$ are odd integers. If $i_1$ and $i_2$ have different parity, then the second component of $v_a$ and $v_b$ give $\left(\frac{(j_2-j_1+1)(s-2)}{4}+\frac{s}{2}\right)\equiv0\pmod{s}$, which is impossible as $\frac{(j_2-j_1+1)(s-2)}{4}+\frac{s}{2}$ is an odd integer and $0\leq j_1, j_2\leq \frac{s}{2}-1$. Therefore, the ordering of the vertices $v_1,v_2,\ldots,v_{rs}$ is distinct and covers all the vertices of $T_{r,s}$. Similarly as in Lemma \ref{Lemma for ac(T_{r,s}) when r,s=0 (mod 4)}, conditions $(a)-(d)$ holds true. Hence the result follows.\qed

\begin{theorem}\label{Theorem for ac(T_{r,s}) when r=1,s=0 (mod 4) or r=3,s=2 (mod 4)} If $r\equiv 1\pmod 4$ and $s\equiv 0\pmod 4$, or $r\equiv 3\pmod 4$ and $s\equiv 2\pmod 4$, then
$$ac(T_{r,s})=\frac{1}{8}\left(r^2s+rs^2-rs-2r-2s+2\right).$$
\end{theorem}
\noindent
\textbf{\textit{Proof.}} Consider the ordering $v_1,v_2,\ldots,v_{rs}$ of $V\left(T_{r,s}\right)$ as given in Lemma \ref{Lemma for ac(T_{r,s}) when r=1,s=0 (mod 4) or r=3,s=2 (mod 4)}. We define a mapping $g : V(T_{r,s})\rightarrow \mathbb{N}\cup \{0\}$ such that 
$
  \,  g(v_1)=0;\,
    g(v_j)=g(v_{j-1}) \,\text{ for } j=2,4,6,\ldots, rs \text{ and }
    g(v_j)=g(v_{j-1})+\left(\frac{r-1}{2}+\frac{s}{2}\right)-d(v_j,v_{j-1})+1 \,\text{ for } j=3,5,7,\ldots, rs-1.
$ \newline
Similarly as Theorem \ref{Theorem for ac(T_{r,s}) when r,s=0 (mod 4)}, we examine that $g$ is an antipodal coloring. As by the definition of $g$, $\epsilon_{j}=0$ for all $j=2,4,\ldots,rs$, and $\epsilon_{j}=1$ for all $j=3,5,\ldots,rs-1$. Now using Lemma \ref{Lemma 2}, the span of $g$ is given as
$g(v_{rs})= (rs-1)\left(\frac{r-1}{2}+\frac{s}{2}\right)-\frac{rs}{2}\left(\frac{r-1}{2}+\frac{s}{2}\right)-\left(\frac{rs}{2}-1\right)\left(\frac{r+s+3}{4}\right)+\left(\frac{rs}{2}-1\right)
    =\frac{1}{8}\left(r^2s+rs^2-rs-2r-2s+2\right).
$

\noindent
Now, $d(v_j,v_{j+1})=\frac{r-1}{2}+\frac{s}{2}=diam(T_{r,s})$ for all $j=1,3,5,\ldots,rs-3$. Also, $d(v_{j+1},v_{j+2})=\frac{r+s+3}{4}=d(v_j,v_{j+2})+\epsilon_{j+1}+\epsilon_{j+2}$ for all $j=1,3,5,\ldots,rs-3$ and $d(v_{rs-1},v_{rs})=\frac{r-1}{2}+\frac{s}{2}=diam(T_{r,s})$. Therefore, by Theorem \ref{Theorem 1.}, $g$ is a minimal antipodal coloring, and hence the result follows.\qed

\begin{lemma} \label{Lemma for ac(T_{r,s}) when r=2,s=0 (mod 4)}
Let $r\equiv 2\pmod 4$ and $s\equiv 0\pmod 4$. Then there exists an ordering $v_1,v_2,\ldots,v_{rs}$ of $V\left(T_{r,s}\right)$ with the following properties
\begin{enumerate}
    \item [$(a)$] $d(v_{2j},v_{2j-1})=\frac{r}{2}+\frac{s}{2}=diam(T_{r,s})$\quad for $j=1,2,\ldots,\frac{rs}{2}.$
    \item [$(b)$] $d(v_{2j+1},v_{2j})=\frac{r+2}{4}+\frac{s}{4}$\quad for  $j=1,2,\ldots,\frac{rs}{2}-1.$
    \item [$(c)$] $d(v_{j},v_{j-2})= \frac{r-2}{4}+\frac{s}{4}$\quad for $j=3,4,\ldots,rs.$ 
    \item [$(d)$] $d(v_j,v_{j-3})\begin{cases} =\frac{r+2}{4}+\frac{s}{4} & \text{ if } j=4,6,8,\ldots,rs, \\
    \geq1 & \text{ if } j=5,7,\ldots,rs-1.
    \end{cases}$
\end{enumerate}
\end{lemma}
\noindent
\textbf{\textit{Proof.}} For each $j\in \{0,1,\ldots, s-1\}$, let $V_j=\{(i,j) : 0\leq i\leq r-1\}$. Then clearly $\{V_0, V_1,\ldots, V_{s-1}\}$ is a partition of $V(T_{r,s})$. Now we define an ordering of $V\left(T_{r,s}\right)$ as given below. For $0\leq i\leq \frac{r}{2}-1$, $v_{4i+1}=\left(\frac{i(r-2)}{2},0\right),\,
        v_{4i+2}=\left(\frac{i(r-2)}{2}+\frac{r}{2},\frac{s}{2}\right),\,
        v_{4i+3}=\left(\frac{i(r-2)}{2}+\frac{r-2}{4},\frac{s}{4}\right),\,
        v_{4i+4}=\left(\frac{i(r-2)}{2}+\frac{3r-2}{4},\frac{3s}{4}\right),\,\\
        v_{2r+4i+1}=\left(\frac{i(r+2)}{2},\frac{s}{2}\right),\,
        v_{2r+4i+2}=\left(\frac{i(r+2)}{2}+\frac{r}{2},0\right),\,
        v_{2r+4i+3}=\left(\frac{i(r+2)}{2}+\frac{3r+2}{4},\frac{3s}{4}\right),\,
        v_{2r+4i+4}=\left(\frac{i(r+2)}{2}+\frac{r+2}{4},\frac{s}{4}\right)$, and $v_{4rj+l+1}=v_{l+1}+(j,j)$ for $0\leq l\leq 4r-1 \text{ and }1\leq j\leq \frac{s}{4}-1. $
Now we will show that the above ordering of $V\left(T_{r,s}\right)$ are distinct, and they cover all the vertices of $T_{r,s}$. Since $\gcd\left(\frac{r-2}{2},r\right)=2$, by Lemma \ref{Lemma 4.1.1}, the elements of the set $\{v_{4i+1} : 0\leq i\leq \frac{r}{2}-1\}$ are distinct. Furthermore, applying Lemma \ref{Lemma 4.1.1}, it can be deduced that for every $j\in \{2,3,4\}$, all the elements in the set $\{v_{4i+j} : 0\leq i\leq \frac{r}{2}-1\}$ are distinct. Also, for every $j\in\{1,2,3,4\}$, $\{v_{4i+j} : 0\leq i\leq \frac{r}{2}-1\}\subset V_{\frac{ls}{4}}$ for some $l\in\{1,2,3,4\}$. So $\bigcup\limits_{j=1}^4\{v_{4i+j} : 0\leq i\leq \frac{r}{2}-1\}\subset\bigcup\limits_{j=1}^4V_{\frac{js}{4}}$. Hence, $v_{a}\neq v_{b}$ for any distinct integers $a,b\in\{1,2,\ldots,2r\}$. Similarly, we check that $v_{2r+i}\neq v_{2r+j}$, for $1\leq i,j\leq 2r$ and $i\neq j$. We note that for $i\neq j$ and $i,j\in\{1,2,\ldots,2r\}$, if $v{i}\neq v_{2r+j}$ then at least one of the following is true,
    \begin{itemize}
        \item [$(i)$] $\frac{i(r-2)}{2} \pmod r=\left(\frac{j(r+2)}{2}+\frac{r}{2}\right) \pmod r.$
        \item [$(ii)$] $\left(\frac{i(r-2)}{2}+\frac{r}{2}\right) \pmod r=\frac{j(r+2)}{2} \pmod r.$
        \item [$(iii)$] $\left(\frac{i(r-2)}{2}+\frac{r-2}{4}\right) \pmod r=\left(\frac{j(r+2)}{2}+\frac{r+2}{4}\right) \pmod r.$
        \item [$(iv)$] $\left(\frac{i(r-2)}{2}+\frac{r-2}{4}+\frac{r}{2}\right) \pmod r=\left(\frac{j(r+2)}{2}+\frac{3r+2}{4}\right) \pmod r.$
    \end{itemize}
    As $\frac{i(r-2)}{2} \pmod r$ is an even integer and $\left(\frac{j(r+2)}{2}+\frac{r}{2}\right) \pmod r$ is an odd integer, so $(i)$ is not true. For the same reason, $(ii)$ is false. Also, if $r\equiv 2 \pmod 8$ or $r\equiv 6 \pmod 8$ then $\left(\frac{i(r-2)}{2}+\frac{r-2}{4}\right) \pmod r$ is even or odd, and $\left(\frac{j(r+2)}{2}+\frac{r+2}{4}\right) \pmod r$ is odd or even, respectively. So $(iii)$ is also not true. Similarly, $(iv)$ is not true. Therefore $v_{i}\neq v_{j}$, for $i\neq j$ and $i,j\in\{1,2,\ldots,4r\}$. Now, clearly $\{v_{i} : 1\leq i\leq 4r\}=\bigcup\limits_{j=1}^4V_{\frac{js}{4}}$. Moreover, we observe that $\{v_{4rj+l+1} : 0\leq l\leq 4r-1\}=\{v_{l+1}+(j,j) : 0\leq l\leq 4r-1\}=V_j\bigcup V_{\frac{s}{4}+j}\bigcup V_{\frac{s}{2}+j}\bigcup V_{\frac{3s}{4}+j}$ for each $j=0,1,\ldots \frac{s}{4}-1$. Therefore the ordering of the vertices $v_1,v_2,\ldots,v_{rs}$ are distinct, and they cover all the vertices of $T_{r,s}$. Similarly as Lemma \ref{Lemma for ac(T_{r,s}) when r,s=0 (mod 4)},  the conditions $(a)-(d)$ holds true. Hence, the result follows.\qed

\begin{theorem}\label{Theorem for ac(T_{r,s}) when r=2,s=0 (mod 4)} If $r\equiv 2\pmod 4$ and $s\equiv 0\pmod 4$, then
$$ac(T_{r,s})=\frac{1}{8}\left(r^2s+rs^2+2rs-2r-2s-8\right).$$
\end{theorem}
\noindent
\textbf{\textit{Proof.}} Consider the ordering $v_1,v_2,\ldots,v_{rs}$ of $V\left(T_{r,s}\right)$ as given in the Lemma \ref{Lemma for ac(T_{r,s}) when r=3,s=0 (mod 4)}. We define a mapping $g : V(T_{r,s})\rightarrow \mathbb{N}\cup \{0\}$ such that 
$
  \,  g(v_1)=0;
    g(v_j)=g(v_{j-1}) \,\text{ for } j=2,4,6,\ldots, rs \text{ and }
    g(v_j)=g(v_{j-1})+\left(\frac{r}{2}+\frac{s}{2}\right)-d(v_j,v_{j-1})+1 \,\text{ for } j=3,5,7,\ldots,rs-1.
$\newline
Similarly as Theorem \ref{Theorem for ac(T_{r,s}) when r,s=0 (mod 4)}, we examine that $g$ is an antipodal coloring. As by the definition of $g$, $\epsilon_{j}=0$ for all $j=2,4,\ldots,rs$, and $\epsilon_{j}=1$ for all $j=3,5,\ldots,rs-1$. Now using Lemma \ref{Lemma 2}, the span of $g$ is given as
$ g(v_{rs})=(rs-1)\left(\frac{r}{2}+\frac{s}{2}\right)-\frac{rs}{2}\left(\frac{r}{2}+\frac{s}{2}\right)-\left(\frac{rs}{2}-1\right)\left(\frac{r+2}{4}+\frac{s}{4}\right)+\left(\frac{rs}{2}-1\right)\cdot 1
    =\frac{1}{8}\left(r^2s+rs^2+2rs-2r-2s-8\right).
$

\noindent
Now $d(v_j,v_{j+1})=\frac{r}{2}+\frac{s}{2}=diam(T_{r,s})$ for all $j=1,3,5,\ldots,rs-3$. Also, $d(v_{j+1},v_{j+2})=\frac{r+2}{4}+\frac{s}{4}=d(v_j,v_{j+2})+\epsilon_{j+1}+\epsilon_{j+2}$ for all $j=1,3,5,\ldots,rs-3$, and $d(v_{rs-1},v_{rs})=\frac{r}{2}+\frac{s}{2}=diam(T_{r,s})$. Therefore, by Theorem \ref{Theorem 1.}, $g$ is a minimal antipodal coloring. Hence, the result follows.\qed

\begin{lemma} \label{Lemma for ac(T_{r,s}) when r=3,s=0 (mod 4)}
Let  $r\equiv 3\pmod 4$ and $s\equiv 0\pmod 4$. Then there exists an ordering $v_1,v_2,\ldots,v_{rs}$ of $V\left(T_{r,s}\right)$ with the following properties
\begin{enumerate}
    \item [$(a)$] $d(v_{2j},v_{2j-1})=\frac{r-1}{2}+\frac{s}{2}=diam(T_{r,s})$\quad for $j=1,2,\ldots,\frac{rs}{2}.$
    \item [$(b)$] $d(v_{2j+1},v_{2j})=\frac{r+1}{4}+\frac{s}{4}$\quad for  $j=1,2,\ldots,\frac{rs}{2}-1.$
    \item [$(c)$] $d(v_{j},v_{j-2})= \begin{cases} \frac{r-3}{4}+\frac{s}{4} & \text{ if } j=3,6,7,10,\ldots,rs-2,rs-1, \\
    \frac{r+1}{4}+\frac{s}{4} & \text{ if } j=4,5,8,9,\ldots,rs-3,rs.
    \end{cases}$ 
    \item [$(d)$] $d(v_j,v_{j-3})\begin{cases} =\frac{r+1}{4}+\frac{s}{4} & \text{ if } j=4,6,8,\ldots,rs, \\
    \geq\frac{s}{2}-1 & \text{ if } j=5,7,\ldots,rs-1.
    \end{cases}$
\end{enumerate}
\end{lemma}
\noindent
\textbf{\textit{Proof.}} Since $r\equiv 3\pmod 4$ and $s\equiv 0\pmod 4$, so consider an ordering of $V\left(T_{r,s}\right)$ as:
    $
        \,v_{4i+1}=\left(\frac{i(r-1)}{2},0\right),\,
        v_{4i+2}=\left(\frac{i(r-1)}{2}+\frac{r-1}{2},\frac{s}{2}\right),\,
        v_{4i+3}=\left(\frac{i(r-1)}{2}+\frac{r-3}{4},\frac{s}{4}\right),
        v_{4i+4}=\left(\frac{i(r-1)}{2}+\frac{3r-1}{4},\frac{3s}{4}\right)$ for $0\leq i\leq r-1$, and $
        v_{4rj+l+1}=v_{l+1}+(j,j) \text{ for } 0\leq l\leq 4r-1 \text{ and }1\leq j\leq \frac{s}{4}-1.
    $
    Similarly as Lemma \ref{Lemma for ac(T_{r,s}) when r,s=0 (mod 4)}, we can show that the result is true.\qed

\begin{theorem}\label{Theorem for ac(T_{r,s}) when r=3,s=0 (mod 4)} If $r\equiv 3\pmod 4$ and $s\equiv 0\pmod 4$, then
$$ac(T_{r,s})=\frac{1}{8}\left(r^2s+rs^2-rs-2r-2s+6\right).$$
\end{theorem}
\noindent
\textbf{\textit{Proof.}} Consider the ordering $v_1,v_2,\ldots,v_{rs}$ of $V\left(T_{r,s}\right)$ as given in the Lemma \ref{Lemma for ac(T_{r,s}) when r=3,s=0 (mod 4)}. We define a mapping $g : V(T_{r,s})\rightarrow \mathbb{N}\cup \{0\}$ such that 
$
 \,   g(v_1)=0;
    g(v_j)=g(v_{j-1}) \,\text{ for } j=2,4,6,\ldots, rs;\,
    g(v_j)=g(v_{j-1})+\left(\frac{r-1}{2}+\frac{s}{2}\right)-d(v_j,v_{j-1})+1 \,\text{ for } j=3,7,\ldots,rs-1;\,
    g(v_j)=g(v_{j-1})+\left(\frac{r-1}{2}+\frac{s}{2}\right)-d(v_j,v_{j-1}) \,\text{ for } j=5,9,\ldots, rs-3.
$\newline
Similarly as Theorem \ref{Theorem for ac(T_{r,s}) when r,s=0 (mod 4)}, we can examine that $g$ is an antipodal coloring. As by the definition of $g$, $\epsilon_{j}=0$ for all $j=2,4,5,6,8,\ldots,rs$, and $\epsilon_{j}=1$ for all $j=3,7,\ldots,rs-3$. Now using Lemma \ref{Lemma 2}, the span of $g$ is given as
$g(v_{rs})=(rs-1)\left(\frac{r-1}{2}+\frac{s}{2}\right)-\frac{rs}{2}\left(\frac{r-1}{2}+\frac{s}{2}\right)-\left(\frac{rs}{2}-1\right)\left(\frac{r+1}{4}+\frac{s}{4}\right)+\left(\frac{rs}{4}\right)\cdot 1
    =\frac{1}{8}\left(r^2s+rs^2-rs-2r-2s+6\right).
$

\noindent
Now $d(v_j,v_{j+1})=\frac{r-1}{2}+\frac{s}{2}=diam(T_{r,s})$ for all $j=1,3,5,\ldots,rs-3$. Also, $d(v_{j+1},v_{j+2})=\frac{r+1}{4}+\frac{s}{4}=d(v_j,v_{j+2})+\epsilon_{j+1}+\epsilon_{j+2}$ for all $j=1,3,5,\ldots,rs-3$, and $d(v_{rs-1},v_{rs})=\frac{r-1}{2}+\frac{s}{2}=diam(T_{r,s})$. Therefore, by Theorem \ref{Theorem 1.}, $g$ is a minimal antipodal coloring. Hence, the result follows.\qed
\begin{lemma}\label{Lemma for ac(T_{r,s}) when r,s=2 (mod 4)}
Let $r\equiv 2\pmod 4$ and $s\equiv 2\pmod 4$. If $r\equiv 2\pmod 8$ and $s\equiv 2\pmod 8$, or $r\equiv 6\pmod 8$ and $s\equiv 6\pmod 8$, then there exists an ordering $v_1,v_2,\ldots,v_{rs}$ of $V\left(T_{r,s}\right)$ with the following properties
\begin{enumerate}
    \item [$(a)$] $d(v_{2j},v_{2j-1})=\frac{r}{2}+\frac{s}{2}=diam(T_{r,s})$\quad for $j=1,2,\ldots,\frac{rs}{2}.$
    \item [$(b)$] $d(v_{2j+1},v_{2j})=\frac{r+2}{4}+\frac{s-2}{4}$\quad for  $j=1,2,\ldots,\frac{rs}{2}-1.$
    \item [$(c)$] $d(v_{j},v_{j-2})= \frac{r-2}{4}+\frac{s+2}{4}$\quad for $j=3,4,\ldots,rs.$ 
    \item [$(d)$] $d(v_j,v_{j-3})\begin{cases} =\frac{r+2}{4}+\frac{s-2}{4} & \text{ if } j=4,6,8,\ldots,rs, \\
    \geq1 & \text{ if } j=5,7,\ldots,rs-1.
    \end{cases}$
\end{enumerate}
If $r\equiv 6\pmod 8$ and $s\equiv 2\pmod 8$, then there exists an ordering $v_1,v_2,\ldots,v_{rs}$ of $V\left(T_{r,s}\right)$ with the following properties
\begin{enumerate}
    \item [$(a)$] $d(v_{2j},v_{2j-1})=\frac{r}{2}+\frac{s}{2}=diam(T_{r,s})$\quad for $j=1,2,\ldots,\frac{rs}{2}.$
    \item [$(b)$] $d(v_{2j+1},v_{2j})=\begin{cases}\frac{r+2}{4}+\frac{s-2}{4} & \text{ if } j=1,2,\ldots,\frac{rs}{2}-1, \\ & \text{ except }j=\frac{s}{2},s,\frac{3s}{2},\ldots,\frac{rs}{2}-\frac{s}{2},\\
    \frac{r+2}{4}+\frac{s+2}{4} & \text{ if } j=\frac{s}{2},s,\frac{3s}{2},\ldots,\frac{rs}{2}-\frac{s}{2}.
    \end{cases}$
    \item [$(c)$] $d(v_j,v_{j-2})=\begin{cases} \frac{r-2}{4}+\frac{s+2}{4} & \text{ if } j=3,4,\ldots,rs,  \\ & \text{except }j=s+1,s+2,2s+1,\ldots,rs-s+2,\\
    \frac{r-2}{4}+\frac{s-2}{4} & \text{if } j=s+1,s+2,2s+1,2s+2,\ldots,rs-s+2.
    \end{cases}$
    \item [$(d)$] $d(v_j,v_{j-3})\begin{cases}= \frac{r+2}{4}+\frac{s-2}{4} & \text{ if } j=4,8,\ldots,rs, \\
    \geq 1 & \text{ if } j=5,7,\ldots,rs-1.
    \end{cases}$
\end{enumerate}
\end{lemma}
\noindent
\textbf{\textit{Proof.}} Since $T_{r,s}$ is isomorphic to $T_{s,r}$, so if the result holds for $r\equiv 6\pmod{8}$ and $s\equiv 2\pmod{8}$ then it holds for $r\equiv 2\pmod{8}$ and $s\equiv 6\pmod{8}$ also. Now we divide the proof of this lemma into two cases as follows:\\

\noindent
    \textbf{Case I:} $r\equiv 6\pmod 8$ and $s\equiv 2\pmod 4$.\\
    Consider an ordering of $V\left(T_{r,s}\right)$ as \newline
   $
        v_{sj+2i+1}=\begin{cases}\left(\frac{j(r-2)}{4},\frac{i(s+2)}{4}\right) & \text{ if } i\equiv0\pmod2, \\
    \left(\frac{(j+1)(r-2)}{4},\frac{i(s+2)}{4}\right) & \text{ if } i\equiv1\pmod2
    \end{cases} \newline
        v_{sj+2i+2}=\begin{cases}\left(\frac{j(r-2)}{4}+\frac{r}{2},\frac{i(s+2)}{4}+\frac{s}{2}\right) & \text{ if } i\equiv0\pmod2, \\
    \left(\frac{(j+1)(r-2)}{4}+\frac{r}{2},\frac{i(s+2)}{4}+\frac{s}{2}\right) & \text{ if } i\equiv1\pmod2
    \end{cases}
   $\newline
    for $i=0,1,\ldots,\frac{s}{2}-1$ and $j=0,1,\ldots,r-1$. Then, we continue similarly as Case II of Lemma \ref{Lemma for ac(T_{r,s}) when r=1,s=0 (mod 4) or r=3,s=2 (mod 4)}, and hence the result follows.\\
   
   \noindent
    \textbf{Case II:} $r\equiv 2\pmod 8$ and $s\equiv 2\pmod 8$.\\
    Consider an ordering of $V\left(T_{r,s}\right)$ as \newline
  $
        v_{sj+2i+1}=\begin{cases}\left(\frac{j(r+2)}{4},\frac{i(s+2)}{4}\right) & \text{ if } i\equiv0\pmod2, \\
    \left(\frac{j(r+2)}{4}+\frac{r-2}{4},\frac{i(s+2)}{4}\right) & \text{ if } i\equiv1\pmod2
    \end{cases}\newline
        v_{sj+2i+2}=\begin{cases}\left(\frac{j(r+2)}{4}+\frac{r}{2},\frac{i(s+2)}{4}+\frac{s}{2}\right) & \text{ if } i\equiv0\pmod2, \\
    \left(\frac{j(r+2)}{4}+\frac{r-2}{4}+\frac{r}{2},\frac{i(s+2)}{4}+\frac{s}{2}\right) & \text{ if } i\equiv1\pmod2
    \end{cases}
 $ \\
    for $i=0,1,\ldots,\frac{s}{2}-1$ and $j=0,1,\ldots,r-1$. Then, we continue similarly as Case II of Lemma \ref{Lemma for ac(T_{r,s}) when r=1,s=0 (mod 4) or r=3,s=2 (mod 4)}, and hence the result follows.\qed
\begin{theorem}\label{Theorem for ac(T_{r,s}) when r=2,s=2 (mod 4)} If $r\equiv 2\pmod 4$ and $s\equiv 2\pmod 4$, then \\
$
    ac(T_{r,s})=\begin{cases}
\frac{1}{8}\left(r^2s+rs^2-2r-2s\right)& \text{if } r\equiv 2\pmod 8\text{ and }s\equiv 2\pmod 8,\\&\text{or }r\equiv 6\pmod 8\text{ and }s\equiv 6\pmod 8,\\
\frac{1}{8}\left(r^2s+rs^2+6r-2s-8\right)& \text{ if }r\equiv 6\pmod 8 \text{ and } s\equiv 2\pmod 8.\end{cases}
$
\end{theorem}
\noindent
\textbf{\textit{Proof.}} Consider the ordering $v_1,v_2,\ldots,v_{rs}$ of $V\left(T_{r,s}\right)$ as given in the Lemma \ref{Lemma for ac(T_{r,s}) when r,s=2 (mod 4)}.\\

\noindent
\textbf{Case I:} $r\equiv 2\pmod 8 \text{ and } s\equiv 2\pmod 8, \text{ or } r\equiv 6\pmod 8 \text{ and } s\equiv 6\pmod 8$.\\
Define a mapping $g : V(T_{r,s})\rightarrow \mathbb{N}\cup \{0\}$ such that 
$
   \, g(v_1)=0;\,
    g(v_j)=g(v_{j-1}) \,\text{ for } j=2,4,6,\ldots, rs;\,\\
    g(v_j)=g(v_{j-1})+\left(\frac{r}{2}+\frac{s}{2}\right)-d(v_j,v_{j-1}) \,\text{ for } j=3,5,7,\ldots,rs-1.
$
Similarly as Theorem \ref{Theorem for ac(T_{r,s}) when r,s=0 (mod 4)}, we can examine that $g$ is an antipodal coloring. As by definition of $g$, $\epsilon_{j}=0$ for all $j=2,3,4,\ldots,rs$. Now using Lemma \ref{Lemma 2}, the span of $g$ is given as
$g(v_{rs})=(rs-1)\left(\frac{r}{2}+\frac{s}{2}\right)-\frac{rs}{2}\left(\frac{r}{2}+\frac{s}{2}\right)-\left(\frac{rs}{2}-1\right)\left(\frac{r+2}{4}+\frac{s-2}{4}\right)
    =\frac{1}{8}\left(r^2s+rs^2-2r-2s\right).
$\\

\noindent
\textbf{Case II:} $r\equiv 6\pmod 8 \text{ and } s\equiv 2\pmod 8$.\\
Define a mapping $g : V(T_{r,s})\rightarrow \mathbb{N}\cup \{0\}$ such that 
$
  \,  g(v_1)=0;\,
    g(v_j)=g(v_{j-1}) \,\text{ for } j=2,4,6,\ldots, rs;\,
    g(v_j)=g(v_{j-1})+\left(\frac{r}{2}+\frac{s}{2}\right)-d(v_j,v_{j-1}) \,\,\text{ for } j=3,5,7,\ldots,rs-1  \text{ except } j=s+1,2s+1,\ldots,rs-s+1;\,
    g(v_j)=g(v_{j-1})+\left(\frac{r}{2}+\frac{s}{2}\right)-d(v_j,v_{j-1})+2 \,\text{ for } j=s+1,2s+1,\ldots,rs-s+1.
$
Similarly as Theorem \ref{Theorem for ac(T_{r,s}) when r,s=0 (mod 4)}, we can examine that $g$ is an antipodal coloring. As by the definition of $g$, $\epsilon_{j}=0$ for all $j=2,3,4,\ldots,rs$ except $j=s+1,2s+1,\ldots,rs-s+1$, and $\epsilon_{j}=1$ for all $j=s+1,2s+1,\ldots,rs-s+1$. Now using Lemma \ref{Lemma 2}, the span of $g$ is given as
$ g(v_{rs})=(rs-1)\left(\frac{r}{2}+\frac{s}{2}\right)-\frac{rs}{2}\left(\frac{r}{2}+\frac{s}{2}\right)-\left(\frac{rs}{2}-1-r+1\right)\left(\frac{r+2}{4}+\frac{s-2}{4}\right)+(r-1)\left(\frac{r+2}{4}+\frac{s+2}{4}\right)+2(r-1)
    =\frac{1}{8}\left(r^2s+rs^2+6r-2s-8\right).
$ 

\noindent
Now for both the cases, $d(v_j,v_{j+1})=\frac{r}{2}+\frac{s}{2}=diam(T_{r,s})$ for all $j=1,3,5,\ldots,rs-3$. Also, $d(v_{j+1},v_{j+2})=d(v_j,v_{j+2})+\epsilon_{j+1}+\epsilon_{j+2}$ for all $j=1,3,5,\ldots,rs-3$, and $d(v_{rs-1},v_{rs})=\frac{r}{2}+\frac{s}{2}=diam(T_{r,s})$. Therefore, by Theorem \ref{Theorem 1.}, $g$ is a minimal antipodal coloring in both cases. Hence, the result follows.\qed

\begin{lemma} \label{Lemma for ac(T_{r,s}) when r=1,s=2 (mod 4)}
Let $r\equiv 1\pmod 4$ and $s\equiv 2\pmod 4$. Then there exists an ordering $v_1,v_2,\ldots,v_{rs}$ of $V\left(T_{r,s}\right)$ with the following properties
\begin{enumerate}
    \item [$(a)$] $d(v_{2j},v_{2j-1})=\frac{r-1}{2}+\frac{s}{2}=diam(T_{r,s})$\quad for $j=1,2,\ldots,\frac{rs}{2}.$
    \item [$(b)$] $d(v_{2j+1},v_{2j})=\frac{r+3}{4}+\frac{s+2}{4}$\quad for  $j=1,2,\ldots,\frac{rs}{2}-1.$
    \item [$(c)$] $d(v_{j},v_{j-2})= \frac{r-1}{4}+\frac{s-2}{4}$\quad for $j=3,4,\ldots,rs.$ 
    \item [$(d)$] $d(v_j,v_{j-3})\begin{cases}= \frac{r-1}{4}+\frac{s+2}{4} & \text{ if } j=4,6,8,\ldots,rs, \\
    \geq 2 & \text{ if } j=5,7,\ldots,rs-1.
    \end{cases}$
\end{enumerate}
\end{lemma}
\noindent
\textbf{\textit{Proof.}} Consider an ordering of the $V\left(T_{r,s}\right)$ as\newline
    $
        v_{2rj+2i+1}=\begin{cases}\left(\frac{i(r-1)}{4},\frac{j(s-2)}{4}\right) & \text{ if } i\equiv0\pmod2, \\
    \left(\frac{(i(r-1)}{4},\frac{(j+1)(s-2)}{4}\right) & \text{ if } i\equiv1\pmod2
    \end{cases}\newline
        v_{2rj+2i+2}=\begin{cases}\left(\frac{i(r-1)}{4}+\frac{r+1}{2},\frac{j(s-2)}{4}+\frac{s}{2}\right) & \text{ if } i\equiv0\pmod2, \\
    \left(\frac{(i(r-1)}{4}+\frac{r+1}{2},\frac{(j+1)(s-2)}{4}+\frac{s}{2}\right) & \text{ if } i\equiv1\pmod2
    \end{cases}
    $ \newline
    for $i=0,1,\ldots,r-1$ and $j=0,1,\ldots,\frac{s}{2}-1$. Then, we continue similarly as Case II of Lemma \ref{Lemma for ac(T_{r,s}) when r=1,s=0 (mod 4) or r=3,s=2 (mod 4)}, and hence the result follows.\qed
\begin{theorem}\label{Theorem for ac(T_{r,s}) when r=1,s=2 (mod 4)} If $r\equiv 1\pmod 4$ and $s\equiv 2\pmod 4$, then 
$$ac(T_{r,s})=\frac{1}{8}\left(r^2s+rs^2+rs-2r-2s-2\right).$$
\end{theorem}
\noindent
\textbf{\textit{Proof.}} Consider the ordering $v_1,v_2,\ldots,v_{rs}$ of $V\left(T_{r,s}\right)$ as given in the Lemma \ref{Lemma for ac(T_{r,s}) when r=1,s=2 (mod 4)}. We define a mapping $g : V(T_{r,s})\rightarrow \mathbb{N}\cup \{0\}$ such that 
$
    g(v_1)=0;\,
    g(v_j)=g(v_{j-1}) \,\text{ for } j=2,4,6,\ldots, rs;\,
    g(v_j)=g(v_{j-1})+\left(\frac{r-1}{2}+\frac{s}{2}\right)-d(v_j,v_{j-1})+2 \,\text{ for } j=3,5,7,\ldots, rs-1.
    $
Similarly as Theorem \ref{Theorem for ac(T_{r,s}) when r,s=0 (mod 4)}, we can examine that $g$ is an antipodal coloring. As from definition of $g$, $\epsilon_{j}=0$ for all $j=2,4,\ldots,rs$, and $\epsilon_{j}=2$ for all $j=3,5,\ldots,rs-1$. Now using Lemma \ref{Lemma 2}, the span of $g$ is given as
$g(v_{rs})=(rs-1)\left(\frac{r-1}{2}+\frac{s}{2}\right)-\frac{rs}{2}\left(\frac{r-1}{2}+\frac{s}{2}\right)-\left(\frac{rs}{2}-1\right)\left(\frac{r+3}{4}+\frac{s+2}{4}\right)+2\left(\frac{rs}{2}-1\right)
    =\frac{1}{8}\left(r^2s+rs^2+rs-2r-2s-2\right).
$ 

\noindent
Now $d(v_j,v_{j+1})=\frac{r-1}{2}+\frac{s}{2}=diam(T_{r,s})$ for all $j=1,3,5,\ldots,rs-3$. Also, $d(v_{j+1},v_{j+2})=\frac{r+3}{4}+\frac{s+2}{4}=d(v_j,v_{j+2})+\epsilon_{j+1}+\epsilon_{j+2}$ for all $j=1,3,5,\ldots,rs-3$, and $d(v_{rs-1},v_{rs})=\frac{r-1}{2}+\frac{s}{2}=diam(T_{r,s})$. Therefore, by Theorem \ref{Theorem 1.}, $g$ is a minimal antipodal coloring. Hence, the result follows.\qed

\noindent
Next, we determine the lower bound for $ac(T_{r,s})$ when $rs$ is an odd integer. To accomplish this, we require the lemma given below.

\begin{lemma}\label{triameter of T_r,s}\cite{saha2013radio}
    In $T_{r,s}$, $d(u,v)+d(v,w)+d(w,u)\leq r+s$ for any vertices $u,v$ and $w$ of $T_{r,s}$.
\end{lemma}

\begin{theorem}\label{lower bound when both r and s are odd} If $r\geq3$ and $s\geq3$ are both odd integers, then
$ac(T_{r,s})\geq \left\lceil\frac{r+s-6}{4}\right\rceil\left(\frac{rs-1}{2}\right)$.    
\end{theorem}
\noindent
\textbf{\textit{Proof.}} Let $g$ be an antipodal coloring of $T_{r,s}$ and $v_1,v_2,\ldots,v_{rs}$ be the ordering of $V\left(T_{r,s}\right)$ in such a way that $g(v_{j+1})>g(v_j)$ for all $j=1,2,\ldots,rs-1$. Then, clearly $g(v_1)=0$ and the span of $g$ is equal to $g(v_{rs})$. So by the definition of antipodal coloring of $g$ of $T_{r,s}$ and $g(v_{j+1})>g(v_j)$ for all $j=1,2,\ldots,rs-1$, we have
\begin{align} 
    g(v_{j+1})-g(v_{j})&\geq diam(T_{r,s})-d(v_{j+1},v_{j}),\label{eq10}\\
    g(v_{j+2})-g(v_{j+1})&\geq diam(T_{r,s})-d(v_{j+2},v_{j+1}),\label{eq11}\\
    g(v_{j+2})-g(v_{j})&\geq diam(T_{r,s})-d(v_{j+2},v_{j}).\label{eq12}
\end{align}
Now adding equations (\ref{eq10}), (\ref{eq11}), (\ref{eq12}), we get
\begin{align} 
    \hspace{-0.3cm}2[g(v_{j+1})-g(v_{j})]&\geq 3 diam(T_{r,s})-[d(v_{j+1},v_{j})+d(v_{j+2},v_{j+1})+d(v_{j+2},v_{j})].\label{eq13}
\end{align}
Now using the Lemma \ref{triameter of T_r,s}, the inequality (\ref{eq13}) becomes
\begin{align}
     2[g(v_{j+1})-g(v_{j})]&\geq 3 diam(T_{r,s})-r-s.\label{eq14}
\end{align}
Since $g(v_{j+1})$ and $g(v_{j})$ are integers, so the inequality (\ref{eq14}) becomes
\begin{align}
    g(v_{j+1})-g(v_{j})&\geq\left\lceil \frac{3diam(T_{r,s})-r-s}{2}\right\rceil \text{ for all }j=1,2,\ldots,rs-2.\label{eq15}
\end{align}
From Lemma \ref{Lemma 2.1.01}, $diam(T_{r,s})=\frac{r-1}{2}+\frac{s-1}{2}$ as $r$ and $s$ are odd integers. Now summing up the inequality for $j=1,2,\ldots,rs-2$, we get
\begin{align*}
    g(v_{rs})-g(v_{1})&\geq\left\lceil \frac{r+s-6}{4}\right\rceil\left(\frac{rs-1}{2}\right).
\end{align*}
Also, $g(v_1)=0$, then we get $g(v_{rs})\geq\left\lceil \frac{r+s-6}{4}\right\rceil\left(\frac{rs-1}{2}\right)$. Hence, the result follows.\qed\\\\
\textbf{Acknowledgment:}
The authors are grateful to the Indian Institute of Technology
Kharagpur, India, for financial support.\\
\noindent
\textbf{Disclosure of Interest:}
The authors report there are no competing interests to declare.



\end{document}